# THE MULTIFRACTAL SPECTRUM OF BROWNIAN INTERSECTION LOCAL TIMES[1]

By Achim Klenke and Peter Mörters

*Johannes Gutenberg-Universität Mainz and University of Bath*

Let $\ell$ be the projected intersection local time of two independent Brownian paths in $\mathbb{R}^d$ for $d = 2, 3$. We determine the lower tail of the random variable $\ell(\mathbb{U})$, where $\mathbb{U}$ is the unit ball. The answer is given in terms of intersection exponents, which are explicitly known in the case of planar Brownian motion. We use this result to obtain the multifractal spectrum, or spectrum of thin points, for the intersection local times.

## 1. Introduction and main results.

1.1. *Aims of the paper.* Intersections of Brownian motion or random walk paths have been studied for quite a long time in probability theory and statistical mechanics. One of the reasons for this interest is that the properties of the intersections are analogous to those of a number of more complicated models in equilibrium statistical physics. There is trivial behavior in all dimensions exceeding a critical dimension, which in our case is $d = 4$, but below the critical dimension there are interesting critical exponents, which determine the universality class of the model and enter into most of its quantitative studies. Rigorous and nonrigorous techniques from mathematical physics, such as renormalization group theory (see, e.g., [1]) and conformal field theory (see, e.g., [6]), have been applied to the model and, more recently, finding the intersection exponents of planar Brownian motion was one of the first problems solved by the rigorous techniques

Received January 2004; revised October 2004.

[1]Supported by the DFG Schwerpunkt programme "Stochastically interacting systems of high complexity," the ESF scientific programme "Random dynamics in spatially extended systems" (RDSES) and Nuffield Foundation Grant NAL/00631/G.

based on the stochastic Löwner evolution devised by Lawler, Schramm and Werner [19, 20, 21].

An interesting geometric characteristic for discrete and continuous models of statistical physics is the *multifractal spectrum*, which originated in the study of turbulence models. Loosely speaking, the multifractal spectrum evaluates the degree of variation in the intensity of a spatial distribution. Calculations of the multifractal spectrum in the physics literature are typically done using a *multifractal formalism*, based on large-deviation heuristics, which emerged in the physics literature in the late 1980s; see, for example, [9]. This formalism allows nonrigorous calculation of the multifractal spectrum in many models, either explicitly or in terms of the critical exponents; see, for example, [5] for a survey from a physicist's point of view.

In some cases the multifractal spectrum could be calculated rigorously. For a precise definition fix a locally finite, fractal measure $\mu$, which may be random or nonrandom. The value $f(a)$ of the multifractal spectrum is the Hausdorff dimension of the set of points $x$ with local dimension

$$(1.1) \qquad \lim_{r \downarrow 0} \frac{\log \mu(B(x,r))}{\log r} = a,$$

where $B(x,r)$ denotes the open ball of radius $r$ centered in $x$. In some cases of interest, the limit in (1.1) has to be replaced by lim inf or lim sup to obtain an interesting nontrivial spectrum. Examples of rigorously verified multifractal spectra for measures arising in probability theory are the occupation measures of stable subordinators, see [10], the states of super-Brownian motion, see [24] and the harmonic measure on a Brownian path, see [15].

The main aim of the present paper is to find the multifractal spectrum of the intersection local time of two independent Brownian paths in $d = 2, 3$. In this example it is not hard to check that the multifractal formalism predicts a *trivial spectrum*, that is, that the set of points in the intersection of the paths where the local dimension differs from the global dimension of the set has dimension zero. This prediction is not correct; it turns out that there is a nontrivial spectrum to the right of the typical value. The spectrum is given in terms of the critical exponents, in this case the intersection exponents. Failure of the multifractal formalism for natural random measures has been observed before; the most notable example is the states of super-Brownian motion, see [24].

Let $W^1, W^2$ be two independent Brownian motions in $\mathbb{R}^d$, $d = 2, 3$, with a joint starting point and running for one unit of time. Let $S = W^1([0,1]) \cap W^2([0,1])$ be the intersection set, which is equipped with a uniform measure, the intersection local time $\ell$. It is well known that, almost surely,

$$\lim_{r \downarrow 0} \frac{\log \ell(B(x,r))}{\log r} = 4 - d \qquad \text{for } \ell\text{-almost every } x \in S,$$

that is, the local dimension of *typical* points equals $4-d$, which is also the Hausdorff dimension of the intersection set $S$. The multifractal spectrum deals with sets of *exceptional* points, of which there may be two types: We call a point $x \in S$ a *thin point* if

$$\limsup_{r \downarrow 0} \frac{\log \ell(B(x,r))}{\log r} > 4-d,$$

noting that this means that, for a sequence $r_n \downarrow 0$, the mass $\ell(B(x,r_n))$ is unusually small, on a logarithmic scale. Analogously, a point would be called a *thick point* if

$$\liminf_{r \downarrow 0} \frac{\log \ell(B(x,r))}{\log r} < 4-d.$$

However, thick points in this sense do not exist; in fact, almost surely,

(1.2) $$\liminf_{r \downarrow 0} \frac{\log \ell(B(x,r))}{\log r} = 4-d \qquad \text{for every } x \in S.$$

The upper bound in (1.2) is easy to show and the lower bound follows from recent work by Dembo, Peres, Rosen and Zeitouni [4] in $d=2$, and by König and Mörters [13] in $d=3$. Indeed, both papers show how a weaker notion of thick points, which operates on a finer scale, can be defined in order to obtain a nontrivial spectrum.

In this paper we are looking at the right end of the multifractal spectrum, asking for the Hausdorff dimension of the set of thin points or, more precisely, the set of points $x$ where, for some $a > 4-d$,

$$\limsup_{r \downarrow 0} \frac{\log \ell(B(x,r))}{\log r} = a.$$

We show that, with $\xi = \xi_d(2,2)$ the intersection exponent, in $d=3$, thin points *do* exist for all values $1 < a \le \xi/(\xi-1)$, but not for any $a > \xi/(\xi-1)$. In $d=2$ we show that thin points exist for all values $2 < a \le 2\xi/(\xi-2)$, but not for any $a > 2\xi/(\xi-2)$. Moreover, we show that the Hausdorff dimension spectrum of thin points is given by

(1.3) $$\dim\left\{x \in S : \limsup_{r \downarrow 0} \frac{\log \ell(B(x,r))}{\log r} = a\right\} = (4-d)\frac{\xi}{a} + 4 - d - \xi,$$

almost surely for all values of $a \ge 4-d$ where the right-hand side is nonnegative. This result is in contrast to the behavior of Brownian occupation measure, the natural analogue for the case $p=1$, where no thin points (in our sense) exist; see [2].

Our paper is also in the tradition of a series of papers by Lawler, who first observed that intersection exponents enter into the Hausdorff dimension of various subsets of the Brownian path. The most famous example is the planar Brownian frontier, which has dimension $2 - \xi_2(2,0)$; others are the set of cut points, which have dimension $2 - \xi_d(1,1)$, and the multifractal spectrum of harmonic measure on a planar Brownian path, which is the concave Legendre transform of the strictly convex function $\lambda \to 2 - \xi_2(2,\lambda)$. See [14] for one of the earliest papers and [17] for a survey.

Before formulating our precise results in Section 1.3 we briefly review the definition and some results about intersection exponents, which are relevant for our work.

1.2. *Intersection exponents.* Suppose $M, N \in \mathbb{N}$ and let $W^1, \ldots, W^{M+N}$ be a family of independent Brownian motions in $\mathbb{R}^d$, $d = 2, 3$, started uniformly on $\partial B(0,1)$. We divide the motions into two packets and look at the union of the paths in each family,

$$\mathfrak{B}^1(R) = \bigcup_{i=1}^{M} W^i([0, \tau_R^i]) \quad \text{and} \quad \mathfrak{B}^2(R) = \bigcup_{i=M+1}^{M+N} W^i([0, \tau_R^i]),$$

where $\tau_R^i$ is the first exit time of $W^i$ from the ball $B(0, R)$.

The event that the two packets of Brownian paths fail to intersect has a decreasing probability as $R \uparrow \infty$. Indeed, it is easy, using subadditivity, to show that there exists a constant $\xi_d(M, N)$ such that

$$(1.4) \qquad \mathbb{P}\{\mathfrak{B}^1(R) \cap \mathfrak{B}^2(R) = \varnothing\} = R^{-\xi_d(M,N) + o(1)} \qquad \text{as } R \uparrow \infty.$$

The numbers $\xi_d(M, N)$ are called the *intersection exponents*. There are natural extensions of the intersection exponents to the case that $M = 0$, the *disconnection exponents*, and to noninteger numbers $M, N > 0$ of Brownian motions, but we do not need this here.

Physicists, for example, [6], have made conjectures about the precise values of the intersection exponents for a long time now. In particular, they conjectured that in the plane many of these exponents are rational numbers. Very recently, Lawler, Schramm and Werner, in a seminal series of papers [19, 20, 21], have been able to verify this rigorously; see [18] for a survey over the complete series of papers. They have shown that

$$(1.5) \qquad \xi_2(M, N) = \frac{(\sqrt{24M + 1} + \sqrt{24N + 1} - 2)^2 - 4}{48}.$$

This gives $\xi_2(2, 2) = 35/12$. As the proof of (1.5) is based on conformal invariance, there is no analogue in $d = 3$. The only value known in dimension $d = 3$ is $\xi_3(1, 2) = \xi_3(2, 1) = 1$. Indeed, there is no reason why $\xi_3(2, 2)$ should be a rational number. The known bounds show that

$$2 = 2\xi_3(2, 1) > \xi_3(2, 2) > \xi_3(2, 1) = 1,$$

where the strict inequalities follow from the strict concavity of $\lambda \mapsto \xi_3(2,\lambda)$ established in [16].

Extensions of the notion of intersection exponents to $p > 2$ packets of Brownian motions in the plane are usually based on the event that there are no *pairwise* intersections of the sets $\mathfrak{B}^1(R), \ldots, \mathfrak{B}^p(R)$. The behavior of intersection local times of $p$ Brownian motions is, of course, related to the event that there are no *joint* intersections of all packets. Again, subadditivity may be used to show that there exists a number $\bar{\xi}_2(M^1, \ldots, M^p) > 0$ such that

$$(1.6) \quad \mathbb{P}\{\mathfrak{B}^1(R) \cap \cdots \cap \mathfrak{B}^p(R) = \varnothing\} = R^{-\bar{\xi}_2(M^1, \ldots, M^p) + o(1)} \qquad \text{as } R \uparrow \infty,$$

where $\mathfrak{B}^1(R), \ldots, \mathfrak{B}^p(R)$ are packets of $M^1, \ldots, M^p$ independent Brownian motions each started uniformly on $\partial B(0,1)$ and running up to the first exit from $B(0,R)$. This class of exponents does not seem to be treated in the literature so far. Note that the exponents $\bar{\xi}$ are smaller than all the exponents discussed before, but it is an open problem to determine their value.

### 1.3. *Main results.*

#### 1.3.1. *Lower tails for intersection local times.*
We now formulate our main results precisely, starting with a result about the *lower tails* of the intersection local times. As we believe that this is of some independent interest, we formulate the result in a somewhat greater generality than needed in our multifractal analysis.

To this end we let $M, N \in \mathbb{N}$ and let $W^1, \ldots, W^{M+N}$ be independent Brownian motions in $\mathbb{R}^d$, $d = 2, 3$, started in the origin. We define the lifetimes of the Brownian motions $W^i$, $1 \leq i \leq M+N$, by

$$\zeta^i := \tau_R^i = \inf\{t > 0 : W^i(t) \notin B(0, R)\},$$

where $R \in (1, \infty]$ may be infinite if $d \geq 3$.

We divide the Brownian motions, as before, in two packets $\mathfrak{B}^1(R)$ and $\mathfrak{B}^2(R)$ of $M$, respectively $N$, motions. On the intersection of the two packets, $S = \mathfrak{B}^1(R) \cap \mathfrak{B}^2(R)$, one can define a natural locally finite measure $\ell$, the (*projected*) *intersection local time*, which can be described symbolically by the formula

$$(1.7) \quad \ell(A) = \sum_{i=1}^{M} \sum_{j=M+1}^{M+N} \int_A dy \int_0^{\zeta^i} ds \int_0^{\zeta^j} dt\, \delta_y(W^i(s)) \delta_y(W^j(t))$$

for $A \subset \mathbb{R}^d$ Borel.

Rigorous constructions of the random measure $\ell$ are reviewed in ([13], Section 2.1). Note that in other sections of this paper we use the same symbol, $\ell$, for intersection local times of Brownian motions running for fixed time. It should always be clear from the context to which situation we are referring.

Let $\mathbb{U} := B(0,1) \subset \mathbb{R}^d$ be the open unit ball in $\mathbb{R}^d$. In [13] the authors determine the upper tails of the random variables $\ell(\mathbb{U})$ in the case $M = N = 1$. It turns out there that $\mathbb{P}\{\ell(\mathbb{U}) > \delta\} \approx \exp(-\theta\delta^{-1/2})$ as $\delta \uparrow \infty$, with the rate $\theta$ given in terms of a variational problem. Our first result shows that the lower tails are fatter, the probabilities $\mathbb{P}\{\ell(\mathbb{U}) < \delta\}$ decaying only polynomially when $\delta \downarrow 0$.

Recall the definition of the intersection exponent $\xi_d(M, N)$ from (1.4).

THEOREM 1.1.
$$\lim_{\delta \downarrow 0} \frac{\log \mathbb{P}\{\ell(\mathbb{U}) < \delta\}}{-\log \delta} = -\frac{\xi_d(M, N)}{4 - d}.$$

An important aspect of this result lies in the fact that the proof also provides an intuitive description of the *strategy* by which the Brownian paths achieve the event $\{\ell(\mathbb{U}) < \delta\}$. Loosely speaking, all Brownian paths run freely until they hit the boundary of the ball $B(0, \delta^{1/(4-d)})$ for the first time. By this time they have accumulated an intersection local time of the order $\delta$. From then on they do not intersect anymore until they exit the unit ball $\mathbb{U}$ for the first time, and after that they never return to the unit ball again. The proof of Theorem 1.1 will be given in Section 2.

1.3.2. *The multifractal spectrum.* We now suppose that $W^1, \ldots, W^p$ are independent Brownian motions in $\mathbb{R}^d$, $d = 2, 3$, started in the origin and running for one unit of time. By classical results of Dvoretzky, Erdős, Kakutani and Taylor (see, e.g., [11, 22] for modern proofs) almost surely the intersection set

(1.8) $$S = W^1([0,1]) \cap \cdots \cap W^p([0,1])$$

contains points different from the origin if and only if $p(d-2) < d$. In these cases the intersection local time $\ell$ is given by the symbolic formula

(1.9) $$\ell(A) = \int_A dy \prod_{i=1}^p \int_0^1 dt \, \delta_y(W^i(t)) \qquad \text{for } A \subset \mathbb{R}^d \text{ Borel.}$$

We focus on the case of two independent Brownian motions $W^1$ and $W^2$ in $\mathbb{R}^d$, $d = 2, 3$, but come back to the case of $p > 2$ motions in the next section.

The fatness of the lower tails observed in Theorem 1.1 is the reason for the existence of thin points, that is, for the fact that there is a nontrivial multifractal spectrum for $\ell$ to the right of the typical value $4 - d$. Our main result determines this spectrum.

DEFINITION 1.2. Denote by

$$\mathcal{T}(a) = \left\{ x \in S : \limsup_{r \downarrow 0} \frac{\log \ell(B(x,r))}{\log r} \geq a \right\},$$

$$\mathcal{T}^s(a) = \left\{ x \in S : \limsup_{r \downarrow 0} \frac{\log \ell(B(x,r))}{\log r} = a \right\}$$

the sets of $a$-thin points, respectively, strictly $a$-thin points in $S$.

Recall from (1.2) that for $d=2$ and $d=3$ there are not thick points; that is, $\mathcal{T}^s(a) = \varnothing$ for $a < 4 - d$ and $\mathcal{T}(a) = \mathcal{T}(4-d)$ for $a \leq 4-d$.

THEOREM 1.3. *Suppose $\ell$ is the intersection local time of two Brownian motions in $\mathbb{R}^d$, for $d = 2, 3$, starting in the origin and running for one unit of time.*

(i) *In $d = 2$ we have*

$$\mathbb{P}\{\mathcal{T}^s(a) \neq \varnothing\} > 0 \quad \text{iff} \quad \mathbb{P}\{\mathcal{T}^s(a) \neq \varnothing\} = 1 \quad \text{iff} \quad 2 \leq a \leq \frac{2\xi_2(2,2)}{\xi_2(2,2) - 2}.$$

*Moreover, for these values of $a$, almost surely,*

$$\dim \mathcal{T}(a) = \dim \mathcal{T}^s(a) = 2\frac{\xi_2(2,2)}{a} + 2 - \xi_2(2,2).$$

(ii) *In $d = 3$ we have*

$$\mathbb{P}\{\mathcal{T}^s(a) \neq \varnothing\} > 0 \quad \text{iff} \quad \mathbb{P}\{\mathcal{T}^s(a) \neq \varnothing\} = 1 \quad \text{iff} \quad 1 \leq a \leq \frac{\xi_3(2,2)}{\xi_3(2,2) - 1}.$$

*Moreover, for these values of $a$, almost surely,*

$$\dim \mathcal{T}(a) = \dim \mathcal{T}^s(a) = \frac{\xi_3(2,2)}{a} + 1 - \xi_3(2,2).$$

The result remains unchanged if the motions are running for any finite amount of time or, in the case $d = 3$, even for infinite time. Note that in the case $d = 2$, by (1.5), we get an explicit multifractal spectrum $f(a) = (1/12)(70/a - 11)$. Let us point out here that the multifractal spectrum for the intersection local times is strictly convex, hence it cannot be found by means of the multifractal formalism, which always predicts concave spectra. The proof of Theorem 1.3 is given in Section 3.

1.3.3. *Intersections of more than two paths.* Recall from [2] that there is no analogous result in the case of a single Brownian path equipped with the occupation measure, as in this case the lower tails are also exponential and thin points fail to exist. There are, however, analogous results for the intersection of any number $p \geq 2$ of Brownian paths in the plane, which we now formulate.

THEOREM 1.4. *Suppose $\ell$ is the intersection local time of $p$ planar Brownian motions, starting in the origin and running for one unit of time. Let $\xi = \bar{\xi}_2(2, \overset{p}{\ldots}, 2) > 0$ be the multiple intersection exponent introduced in* (1.6). *Then,*

$$\mathbb{P}\{\mathcal{T}^s(a) \neq \varnothing\} > 0 \quad \textit{iff} \quad \mathbb{P}\{\mathcal{T}^s(a) \neq \varnothing\} = 1$$

$$\textit{iff} \quad \begin{cases} 2 \leq a \leq 2\xi/(\xi - 2), & \textit{if } \xi > 2, \\ 2 \leq a < \infty, & \textit{if } \xi \leq 2. \end{cases}$$

*Moreover, for these values of $a$,*

$$\dim \mathcal{T}(a) = \dim \mathcal{T}^s(a) = 2\frac{\xi}{a} + 2 - \xi \qquad \textit{almost surely.}$$

For $p = 2$ we have $\bar{\xi}_2(2, 2) = 35/12 > 2$, and hence there is a finite critical value $2\xi/(\xi - 2) = 70/11$ beyond which no $a$-thin points exist. We do *not* know whether the critical value is still finite for larger values of $p$. The proof of Theorem 1.4 is very similar to the proof of Theorem 1.3 and hence details are omitted here.

1.4. *Overview.* We have divided the remainder of this paper into two sections. Section 2 is devoted to the tail asymptotics at zero of the intersection local times, and Section 3 is devoted to the proof of the multifractal spectrum of intersection local time.

In Section 2.1 we show that, if two Brownian paths intersect, they instantly produce some positive amount of intersection local time. This is a nontrivial fact, as times when the paths intersect are *not* stopping times for at least one of the Brownian motions. The exact statement, Proposition 2.3, is a crucial ingredient in the proof of the upper bound in the tail asymptotics. The proof uses the nonpolarity of one Brownian path to show that if *one* intersection occurs, we immediately have a large number of intersections. We then use self-similarity of the paths and a *decoupling technique* to argue that this necessarily leads to positive intersection local time.

In Section 2.2 we give the proof of the upper bound in Theorem 1.1. The proof is based on a *coarse graining* technique. We split the Brownian paths into pieces using suitably defined stopping times. Whenever two pieces intersect, by Proposition 2.3, some positive amount of intersection local time is produced. As the total amount of intersection local time allowed is small, many pieces do not intersect, leading to the upper bound in the probability. This fairly rough argument is only able to give *logarithmic* asymptotics in Theorem 1.1; one would conjecture that $\mathbb{P}\{\ell(\mathbb{U}) < \delta\}$ can be estimated up to constants by a power of $\delta$, but our technique fails to achieve this.

Section 2.3 contains the proof of the lower bound in Theorem 1.1. Here we only have to show that the following strategy (explained already after the theorem) is successful: Put no restrictions on the Brownian motions until they leave the small ball of radius $\delta^{1/(4-d)}$ for the first time, but demand that their paths do not intersect afterwards. A difficulty lies in the fact that the paths might return to this ball and produce more intersection local time by intersecting with "old" pieces of the path. We solve this problem by giving a separate bound for the intersection mass in the small ball. An alternative would be to use results of Lawler [14] to control the probability that Brownian motions do not intersect and simultaneously do not return to the small ball except in a small neighborhood of their respective starting points, but we have opted for the more self-contained argument here.

In Section 2.4 we establish the connection between our tail asymptotics and the multifractal spectrum, providing "local versions" of the tail results in the form needed in the proof of Theorem 1.3; see Lemmas 2.9 and 2.11. Clearly, one can obtain the tail behavior of *small* balls of radius $r > 0$, by Brownian scaling using that, in law, $r^{d-4}\ell_R(B(0,r)) = \ell_{R/r}(B(0,1))$, where the index at $\ell$ indicates the size of the ball where the Brownian motions are stopped. In Lemma 2.11, instead of looking at a pair of Brownian motions started in the same point, we fix a point $x$ different from the motions' starting points, and give the probability that the intersection local time $\ell(B(x,r))$ is smaller than $r^a$, if both Brownian motions are conditioned on hitting $x$. Heuristically, we split the Brownian paths upon first hitting $\partial B(x, r^{a/(4-d)})$. To the incoming paths we apply a *time-reversal* and end up with two pairs of paths which are approximately independent Brownian motions, started uniformly on $\partial B(x, r^{a/(4-d)})$. Now Theorem 1.1 and Brownian scaling yield that the desired probability is of order $r^{\xi_d(2,2)(1-a/(4-d))}$. Lemma 2.11 is used in the proof of the upper bound in Theorem 1.3, and Lemma 2.9 is a variant tailored for use in the lower bound in Theorem 1.3.

In Section 3.1 we verify the upper bounds in Theorem 1.3; see Proposition 3.1. Given Lemma 2.11 these are relatively standard and based on the *first moment method*, that is, on estimating expectations.

In Section 3.2 we explain the setup of the proof of the lower bounds in Theorem 1.3. Our technique uses *percolation limit sets* $\Gamma[\gamma]$ as test sets to determine the Hausdorff dimension of a fractal. More precisely, if a fractal $A$, in our case the set of thin points, hits a test set with a certain parameter $\gamma$ with positive probability, this gives a $\gamma$-dependent lower bound on the dimension of $A$; see Lemma 3.3. The crucial hitting estimate of $\Gamma[\gamma]$ and the set of thin points is formulated as Proposition 3.5, and in Section 3.2 we only show that the upper bounds in Theorem 1.3 follow from this. The remainder of the paper is then devoted to the proof of Proposition 3.5.

In Section 3.3 we show how to overcome the main obstacle in the proof of Proposition 3.5, the long-range dependence. Note that (other than in the thick points problem) long-range dependence is intrinsic in the problem of thin points: If a ball carries very small intersection local time over some time interval, *any* of the two Brownian motions may always return to that ball and create more intersection local time. The main result of this section, Proposition 3.8, shows that in a suitable sense a large number of dyadic cubes are visited only once by both Brownian motions. The proof uses a *two-scale technique* similar to the one used in [2]. On the *coarse* scale we use a dimension argument to ensure that we have enough cubes which are visited by both Brownian motions and retained in the percolation procedure. Within every such big cube we can independently use arguments on the *fine* scale, based on *decoupling* and delicate *second moment estimates*, which show that we have sufficiently many small scale cubes, which are not revisited before the motions leave the big cube. Finally, ensuring that the motions do not revisit many small cubes after leaving the big cube only needs a *first moment* technique. This section is the technically most demanding part of the paper.

In Section 3.4 we complete the proof of Proposition 3.5. Thanks to Proposition 3.8 one can focus on a *localized* notion of thin points, and use the decoupling technique and Lemma 2.9 to ensure the existence of thin points in the percolation limit set. This final part of the proof follows largely the arguments of [13].

**2. Lower tail asymptotics.** Throughout the paper we use the following notation. For any open or closed sets $A_1, A_2, \ldots$ and $i = 1, \ldots, p$ define

(2.1)
$$\tau^i(A_1) := \inf\{t \geq 0 : W^i(t) \in A_1\},$$
$$\tau^i(A_1, \ldots, A_n) := \begin{cases} \inf\{t \geq \tau^i(A_1, \ldots, A_{n-1}) : W^i(t) \in A_n\}, \\ \qquad \text{if } \tau^i(A_1, \ldots, A_{n-1}) < \infty, \\ \infty, \qquad \text{otherwise.} \end{cases}$$

Further for $n \in \mathbb{N}$ and $R_1, \ldots, R_n > 0$ let

(2.2) $$\tau^i_{R_1,\ldots,R_n} := \tau^i(\partial B(0, R_1), \ldots, \partial B(0, R_n))$$

be the hitting time of $\partial B(0, R_n)$ after having hit (in this order) $\partial B(0, R_1), \ldots, \partial B(0, R_{n-1})$.

For the reader's convenience we recall the following well-known lemma for the hitting time of concentric balls. For a single Brownian motion $W$ let $\tau(A) := \inf\{t \geq 0 : W(t) \in A\}$ and $\tau_r := \tau(\partial B(0, r))$ be the first exit time from $B(0, r)$.

LEMMA 2.1. *Let $r_1 \leq r \leq r_2$ and let $W$ be a Brownian motion started in some point in $\partial B(0, r)$. Then*

$$\mathbb{P}\{\tau_{r_1} < \tau_{r_2}\} = f_d(r, r_1, r_2) := \begin{cases} \dfrac{\log(r/r_2)}{\log(r_1/r_2)}, & \text{if } d = 2, \\ \dfrac{(r_2/r) - 1}{(r_2/r_1) - 1}, & \text{if } d = 3, \end{cases}$$

*and*

$$\mathbb{P}\{\tau_{r_1} < \infty\} = \begin{cases} 1, & \text{if } d = 2, \\ r_1/r, & \text{if } d = 3. \end{cases}$$

The proof is standard and can be found in textbooks, for example, in ([7], Chapter 3). From this statement we get the following useful corollary.

COROLLARY 2.2. *Let $\rho \in (0, 1/2)$ and $x \in \mathbb{R}^d$ and assume $r_1 < \rho r < \rho^2 r_2$. Let $D$, $D_1$ and $D_2$ be each either a ball of radius $r$ (resp. $r_1$ or $r_2$) or a box of sidelength $r$ (resp. $r_1$ or $r_2$), centered in $x$. Further let $y \in \partial D$ and $z \in \partial D_2$ and let $W$ be a Brownian motion started in $y$. Then there exists a constant $\tilde{c} := \tilde{c}(\rho) \in (0, \infty)$ depending only on $\rho$ such that*

$$\begin{aligned}
\frac{1}{\tilde{c}} f_d(r, r_1, r_2) &\leq \mathbb{P}\{\tau(\partial D_1) < \tau(\partial D_2)\} \leq \tilde{c} f_d(r, r_1, r_2), \\
\frac{1}{\tilde{c}} f_d(r, r_2, r_1) &\leq \mathbb{P}\{\tau(\partial D_2) < \tau(\partial D_1)\} \leq \tilde{c} f_d(r, r_2, r_1),
\end{aligned} \quad (2.3)$$

*and, in the case where $D_2$ is a ball,*

$$(2.4) \quad \frac{1}{\tilde{c}} f_d(r, r_1, r_2) \leq \mathbb{P}\{\tau(\partial D_1) < \tau(\partial D_2) | W(\tau(\partial D_2)) = z\} \leq \tilde{c} f_d(r, r_1, r_2).$$

PROOF. Without loss of generality we may assume $x = 0$. Note that, in the case where $D_1$ is a box, changing $D_1$ into a ball of radius $r_1/2$ only decreases the probability we want to estimate. Now

$$\frac{f_3(r, r_1/2, r_2)}{f_3(r, r_1, r_2)} = \frac{r_2/r_1 - 1}{2r_2/r_1 - 1} \geq \frac{\rho^{-2} - 1}{2\rho^{-2} - 1} > 0$$

and

$$\frac{f_2(r, r_1/2, r_2)}{f_2(r, r_1, r_2)} = \frac{\log(r_1/r_2)}{\log(r_1/r_2) - \log(2)} \geq \frac{2\log(\rho)}{2\log(\rho) - \log(2)} > 0.$$

On the other hand, changing $D_1$ into a ball of radius $r_1$ increases the probability we want to estimate. By a similar argument we may assume that $D$ is a ball. For the first equation in (2.3) it is sufficient to apply the same argument once more to $D_2$, and the proof of the second equation is analogous.

The proof of (2.4) requires a little more work. We may now assume that $D_1$, $D$ and $D_2$ are balls. Define the open annulus $A_r := \{u \in \mathbb{R}^d : \|u\| \in (3r/4, 3r/2)\}$. Note that $A_r \cap D_1 = \varnothing$ and $\partial D \subset A_r \subset D_2$. Define the random time

$$T := \sup\{t \geq 0 : W(t) \in \partial D \text{ and } W(s) \in A \text{ for all } s \leq t\},$$

which is the starting time of the first excursion off $\partial D$ that leaves $A$. Note that [with $\mathcal{U}_r$ the normalized Lebesgue measure on $\partial B(0, r)$] there is a constant $c > 0$, independent of $y$, $u$ and $r$, such that

$$c^{-1} \leq \frac{\mathbb{P}_y\{W(T) \in du\}}{\mathcal{U}_r(du)} \leq c.$$

Further, for $u \in \partial D$ let $\mathbb{Q}_u$ be the law of a Brownian motion started in $u$ and conditioned to leave $A$ before it returns to $\partial D$ (if it does),

$$\mathbb{Q}_u := \mathcal{L}_u(W | W(t) \notin D \text{ for all } t \in (0, \tau(A^c))).$$

We can now decompose the Brownian motion path into the piece before and into the piece after $T$ to obtain for any measurable set $B \subset \partial D_2$:

$$\mathbb{P}_y\{\tau(\partial D_1) < \tau(\partial D_2) \text{ and } W(\tau(\partial D_2)) \in B\}$$
$$= \int_{\partial D} \mathbb{P}_y\{W(T) \in du\} \mathbb{Q}_u\{\tau(\partial D_1) < \tau(\partial D_2) \text{ and } W(\tau(\partial D_2)) \in B\}$$
$$\geq \frac{1}{c} \int_{\partial D} \mathcal{U}_r(du) \mathbb{Q}_u\{\tau(\partial D_1) < \tau(\partial D_2) \text{ and } W(\tau(\partial D_2)) \in B\}.$$

By rotational symmetry this equals

$$\frac{1}{c} \mathcal{U}_{r_2}(B) \int_{\partial D} \mathcal{U}_r(du) \mathbb{Q}_u\{\tau(\partial D_1) < \tau(\partial D_2)\}$$
$$= \frac{1}{c} \mathcal{U}_{r_2}(B) \int_{\partial D} \mathcal{U}_r(du) \mathbb{P}_u\{\tau(\partial D_1) < \tau(\partial D_2)\}$$
$$= \frac{1}{c} \mathcal{U}_{r_2}(B) f_d(r, r_1, r_2).$$

Analogously we get $\mathbb{P}_y\{\tau(\partial D_1) < \tau(\partial D_2) \text{ and } W(\tau(\partial D_2)) \in B\} \leq c\mathcal{U}_{r_2}(B) f_d(r, r_1, r_2)$. For a constant $c' \in (0, \infty)$ such that $1/c' \leq \mathbb{P}_y\{W(\tau(\partial D_2)) \in du\}/\mathcal{U}_{r_2}(du) \leq c'$, (2.4) holds with $\tilde{c} = c \cdot c'$. □

2.1. *Intersecting paths produce intersection local time.* A basic principle in the proof of Theorem 1.1 is that, whenever two paths intersect, they immediately produce a positive amount of intersection local time. This statement, Proposition 2.3, is proved using a *decoupling technique*, which is also a fundamental tool in the proof of the lower bound for the multifractal spectrum, performed in Section 3.

PROPOSITION 2.3. *Let $W^1$, $W^2$ be two independent Brownian paths with $W^1(0), W^2(0) \in \mathbb{U}$ and $\tau^1 = \tau^1(\mathbb{U}^c)$, $\tau^2 = \tau^2(\mathbb{U}^c)$ be the first exit times from the unit ball. Moreover let $S = W^1([0, \tau^1]) \cap W^2([0, \tau^2])$ be the intersection of the paths, and let $\ell$ be the intersection local time of the paths stopped at time $\tau^1$, respectively, $\tau^2$. Then*

$$\mathbb{P}\{\ell(\mathbb{U}) > 0 | S \neq \varnothing\} = 1.$$

PROOF. First fix the path $W^1$ and let $A := W^1([0, \tau^1])$. Define a stopping time $\sigma = \inf\{t > 0 : W^2(t) \in A\}$ for $W^2$ and recall that, if $\sigma < \tau^2$, the point $W^2(\sigma)$ is regular for the set $A$, which means that

$$\inf\{t > \sigma : W^2(t) \in A\} = \sigma \quad \text{almost surely.}$$

As points are polar, almost surely, there exists $t \in (\sigma, \tau^2)$ such that $W^2(t) \in A \setminus \{W^2(\sigma)\}$. Hence, given $\delta > 0$ we can find a small $\varepsilon > 0$ such that

$$\mathbb{P}\{\inf\{t > \sigma : W^2(t) \in A \setminus B(W^2(\sigma), \varepsilon)\} < \tau^2 | \sigma < \tau^2\} > 1 - \delta.$$

For every integer $M \geq 2$ we can iterate this procedure $M2^d$ times and, averaging over $W^1$ again, we can find for every $\delta > 0$ an $\varepsilon > 0$ such that the event

$$A_M = \{\text{there exist } x_1, \ldots, x_{M2^d} \in S \text{ with } |x_i - x_j| > \varepsilon \, \forall i \neq j, \, |x_i| < 1 - \varepsilon \, \forall i\}$$

satisfies $\mathbb{P}\{A_M | S \neq \varnothing\} > 1 - \delta$. Hence

$$\mathbb{P}\{\ell(\mathbb{U}) = 0 | S \neq \varnothing\} \leq \mathbb{P}\{S \neq \varnothing\}^{-1} \mathbb{P}(\{\ell(\mathbb{U}) = 0\} \cap A_M) + \delta.$$

It therefore suffices to show that there exists an absolute constant $\tilde{\varepsilon} > 0$ such that, for any large $M$,

(2.5) $$\mathbb{P}(\{\ell(\mathbb{U}) = 0\} \cap A_M) \leq 2^d (1 - \tilde{\varepsilon})^M \stackrel{M \to \infty}{\longrightarrow} 0.$$

To verify (2.5) we write $\mathfrak{D}_k$ for the collection of all dyadic cubes $V = \prod_{i=1}^{d} [k_i/2^k, (k_i + 1)/2^k)$, $k_1, \ldots, k_d$ integers. For each such cube $V$ we denote by $B(V)$ the open ball centered in the center of $V$, of radius $2^{-k}$. We denote

$$\mathfrak{D}_k(\mathbb{U}) := \{V \in \mathfrak{D}_k : B(V) \subset \mathbb{U}\}.$$

Furthermore divide $\mathfrak{D}_k(\mathbb{U})$ into $m = 2^d$ subfamilies $\mathfrak{D}_k(\mathbb{U}, 1), \ldots, \mathfrak{D}_k(\mathbb{U}, m)$ such that $B(V) \cap B(V') = \varnothing$ if $V \neq V'$ are in the same subfamily.

Fix $k$ such that $\sqrt{d} \, 2^{-k} < \varepsilon$. For each $j = 1, \ldots, m$ we now introduce a decoupling $\sigma$-field $\mathcal{F}_k(j)$. The idea is to consider the first entrance $\rho(1)$ of a path into one of the cubes $V(1) \in \mathfrak{D}_k(\mathbb{U}, j)$, then its first exit $\sigma(1)$ [after $\rho(1)$] of the ball $B(V(1))$, after $\sigma(1)$ its first entrance $\rho(2)$ into some new box $V(2)$, and so on. $\mathcal{F}_k(j)$ will then use information of the paths *between* the successive times of leaving $B(V(n))$ and entering $V(n+1)$, $n \in \mathbb{N}$.

For the moment we suppress $j$ in the notation. Formally for $i = 1, 2$ we introduce a sequence of (random) sets $(V^i(n) : n = 1, \ldots, \nu^i)$ and stopping times

$$0 = \sigma^i(0) < \rho^i(1) < \sigma^i(1) < \rho^i(2) < \cdots < \sigma^i(\nu^i) < \tau^i < \rho^i(\nu^i + 1),$$

by

$$\begin{aligned}
\rho^i(1) &:= \inf\{\tau^i(V) : V \in \mathfrak{D}_k(\mathbb{U}, j)\}, \\
\tau^i(V^i(n)) &= \rho^i(n) \qquad [\text{this defines } V^i(n)], \\
\sigma^i(n) &= \tau^i(V^i(n), B(V^i(n))^c) \qquad \text{if } \rho^i(n) < \infty, \\
\rho^i(n+1) &= \inf\{\tau^i(V) : V \in \mathfrak{D}_k(\mathbb{U}, j) \setminus \{V^i(1), \ldots, V^i(n)\}\}, \\
\nu^i &:= \max\{n : \rho^i(n) < \tau^i\}.
\end{aligned}$$

Now define

$$\mathcal{F}_k^i(j) := \sigma(W^i(\sigma^i(n) + t), \ t \in [0, \rho^i(n+1) - \sigma^i(n)], \ n = 0, \ldots, \nu^i),$$

and $\mathcal{F}_k(j) := \mathcal{F}_k^1(j) \vee \mathcal{F}_k^2(j)$. Denote

$$\mathfrak{B}_k(j) := \{V \in \mathfrak{D}_k(\mathbb{U}, j) : \tau^i(V) < \tau^i \text{ for all } i = 1, 2\},$$

and note that the events $\{V \in \mathfrak{B}_k(j)\}$, for $V \in \mathfrak{D}_k$, are in $\mathcal{F}_k(j)$. Also observe that

$$A_M \subset \bigcup_{j=1}^m \{\#\mathfrak{B}_k(j) \geq M\},$$

and that

$$\{\#\mathfrak{B}_k(j) \geq M\} \in \mathcal{F}_k(j) \qquad \text{for every } j.$$

It follows easily from the nontriviality of the intersection local times and the boundedness of the density of the harmonic measure, that there exists an absolute constant $\tilde{\varepsilon} > 0$ such that, for all $k \in \mathbb{N}$, $V \in \mathfrak{D}_k$, $x^1, x^2 \in \partial V$, $y^1, y^2 \in \partial B(V)$,

$$\mathbb{P}_{x^1, x^2}\{\ell_V(B(V)) > 0 | W^1(\tau^1(B(V)^c)) = y^1, W^2(\tau^2(B(V)^c)) = y^2\} > \tilde{\varepsilon},$$

where $\mathbb{P}_{x^1, x^2}$ refers to two Brownian motions $W^1, W^2$ with $W^1(0) = x^1$ and $W^2(0) = x^2$, and $\ell_V$ is the intersection local time of the paths $W^1([0, \tau^1(B(V)^c)])$ and $W^2([0, \tau^2(B(V)^c)])$.

Note that, given $\mathcal{F}_k(j)$, the family of random variables $(\ell_V : V \in \mathfrak{B}_n(j))$ are independent. We can now put this information together and get

$$\mathbb{P}(\{\ell(\mathbb{U}) = 0\} \cap A_M)$$

$$\leq \sum_{j=1}^{m} \mathbb{P}\{\ell_V(B(V)) = 0 \text{ for all } V \in \mathfrak{B}_k(j), \ \#\mathfrak{B}_k(j) \geq M\}$$

$$= \sum_{j=1}^{m} \mathbb{E}[\mathbb{P}\{\ell_V(B(V)) = 0 \text{ for all } V \in \mathfrak{B}_k(j)|\mathcal{F}_k(j)\}\mathbb{1}_{\{\#\mathfrak{B}_k(j) \geq M\}}]$$

$$= \sum_{j=1}^{m} \mathbb{E}\left[\prod_{V \in \mathfrak{B}_k(j)} \mathbb{P}\{\ell_V(B(V)) = 0|\mathcal{F}_k(j)\}\mathbb{1}_{\{\#\mathfrak{B}_k(j) \geq M\}}\right] \leq m(1-\tilde{\varepsilon})^M,$$

which is (2.5), and hence the proof is complete. $\square$

2.2. *Proof of Theorem* 1.1, *upper bound.* The idea of the proof is to use a sequence of stopping times to divide each Brownian path into disjoint pieces. Whenever there is an intersection between matching pieces of the two packets, a certain amount of intersection local time is produced. The task is to establish some form of independence between the pieces and estimate the probability of no intersection between matching pieces.

We need three lemmas to prepare the proof of Theorem 1.1. For $r > 0$ denote by $\mathcal{U}_r$ the uniform distribution on $\partial B(0,r)$. For $x \in B(0,r)$ let $m_{r,x}(dy) = \mathbb{P}\{W(\tau_r) \in dy|W(0) = x\}$ be the harmonic measure on $\partial B(0,r)$ for Brownian motion started in $x$. If $x \in \partial B(0,r)$ and $\alpha > 1$, let

$$(2.6) \quad C_\alpha := \sup_{y \in \partial B(0,\alpha r)} \left.\frac{m_{\alpha r,x}(dz)}{\mathcal{U}_{\alpha r}(dz)}\right|_{z=y}, \qquad c_\alpha := \inf_{y \in \partial B(0,\alpha r)} \left.\frac{m_{\alpha r,x}(dz)}{\mathcal{U}_{\alpha r}(dz)}\right|_{z=y}$$

be the maximal and minimal value of the density of $m_{\alpha r,x}$ with respect to the uniform distribution on $\partial B(0,\alpha r)$. Note that, by Brownian scaling and rotational symmetry, both values depend neither on $x$ nor on $r$. Further note that (by the Markov property of $W$), $\alpha \mapsto C_\alpha$ is decreasing and $\alpha \mapsto c_\alpha$ is increasing. For finite measures $\mu$ and $\nu$ we use the ordering $\mu \leq \nu$ iff $\mu(A) \leq \nu(A)$ for all measurable $A$.

LEMMA 2.4. *Let $\mathcal{L}_x$ and $\mathcal{L}_{\mathcal{U}_r}$ denote the laws of Brownian motion $W$ started in $x \in \partial B(0,r)$, respectively, in a point uniformly distributed on $\partial B(0,r)$. Given a Brownian path $W : [0,\infty) \to \mathbb{R}^d$ and $K > 4$ we define*

$$W^{(r)} : [0, \tau_{Kr/2} - \tau_{2r}] \to \mathbb{R}^d, \qquad W^{(r)}(t) = W(\tau_{2r} + t).$$

*With the notation of (2.6) we have for all $y \in \partial B(0, Kr)$ that*

$$\mathcal{L}_x(W^{(r)}|W(\tau_{Kr}) = y) \leq (C_2^2/c_2)\mathcal{L}_{\mathcal{U}_r}(W^{(r)}).$$

PROOF. Fix $B \subset \partial B(0, Kr)$ Borel and a suitable (say bounded continuous) test function $\Psi : C([0,\infty); \mathbb{R}^d) \to \mathbb{R}$. The strong Markov property and three applications of (2.6) yield

$$\mathbb{E}_x[\Psi(W^{(r)})\mathbb{1}_{W(\tau_{Kr})\in B}]$$
$$= \int_{\partial B(0,2r)} m_{2r,x}(dz)\mathbb{E}[\Psi(W^{(r)})\mathbb{1}_{W(\tau_{Kr})\in B}|W(\tau_{2r})=z]$$
$$\leq C_2 \int_{\partial B(0,2r)} \mathcal{U}_{2r}(dz)\mathbb{E}[\Psi(W^{(r)})\mathbb{1}_{W(\tau_{Kr})\in B}|W(\tau_{2r})=z]$$
$$= C_2 \int_{\partial B(0,r)} \mathcal{U}_r(dy) \int_{\partial B(0,2r)} \mathbb{P}_y\{W(\tau_{2r})\in dz\}$$
$$\qquad \times \mathbb{E}[\Psi(W^{(r)})\mathbb{1}_{W(\tau_{Kr})\in B}|W(\tau_{2r})=z]$$
$$= C_2 \mathbb{E}_{\mathcal{U}_r}[\Psi(W^{(r)})\mathbb{1}_{W(\tau_{Kr})\in B}]$$
$$= C_2 \int_{\partial B(0,Kr/2)} \mathbb{E}_{\mathcal{U}_r}[\Psi(W^{(r)})|W(\tau_{Kr/2})=z]\mathbb{P}\{W(\tau_{Kr/2})\in dz\}m_{Kr,z}(B)$$
$$\leq C_2^2 \mathbb{E}_{\mathcal{U}_r}[\Psi(W^{(r)})]\mathcal{U}_{Kr}(B)$$
$$\leq (C_2^2/c_K)\mathbb{E}_{\mathcal{U}_r}[\Psi(W^{(r)})]\mathbb{P}_x\{W(\tau_{Kr})\in B\},$$

from which the result readily follows. □

Fix an arbitrary small $\varepsilon > 0$ and assume that $r \in (0, 1/4)$ is such that

(2.7) $$r \leq (4^{\xi_d(M,N)}(C_2^2/c_2)^{M+N})^{-1/\varepsilon}.$$

For $i = 1, \ldots, M+N$ and any nonnegative integer $k$, let

$$R_k^i := \tau_{2r^{k+1}}^i, \qquad S_k^i := \tau_{4r^{k+1}}^i, \qquad T_k^i := \tau_{r^k}^i.$$

Note that for Brownian motions with $W^i(0) = 0$ for $i = 1, \ldots, N+M$, we have $R_k^i < S_k^i < T_k^i < R_{k-1}^i < \cdots$ for every positive integer $k$. The idea is to consider the Brownian paths in the intervals $[S_k^i, T_k^i]$, $k \geq 1$, only, and to use the remaining intervals for a decoupling of these paths. Hence, we let $L_k$ be the intersection local time of the packets

$$\mathcal{W}_k^1 := \bigcup_{i=1}^M W^i([S_k^i, T_k^i]) \quad \text{and} \quad \mathcal{W}_k^2 := \bigcup_{i=M+1}^{M+N} W^i([S_k^i, T_k^i]).$$

LEMMA 2.5. $\mathbb{P}\{L_0 = 0\} = \mathbb{P}\{\mathcal{W}_0^1 \cap \mathcal{W}_0^2 = \varnothing\}$.

PROOF. This is immediate from Proposition 2.3 and Brownian scaling. □

In addition to $4r < 1$ and (2.7) we assume that $r$ is small enough such that
$$\mathbb{P}\{L_0 = 0\} = \mathbb{P}\{\mathcal{W}_0^1 \cap \mathcal{W}_0^2 = \varnothing\} \leq \tfrac{1}{2}(4r)^{\xi_d(M,N)-\varepsilon},$$
which is possible by Lemma 2.5, the definition (1.4) of $\xi_d(M,N)$ and Brownian scaling. We let
$$\mathcal{F} := \sigma(W^i(R_k^i), i = 1, \ldots, M+N, \ k \geq 0).$$
Given $\mathcal{F}$, the random variables $L_k$ depend only on packets of Brownian motions with disjoint time intervals and fixed initial and exit points, hence by the strong Markov property the sequence $(L_k)_{k \in \mathbb{N}}$ is independent.

LEMMA 2.6. *Let $X_0, X_1, \ldots$ be independent copies of $L_0$ and define $C_1 := (C_2^2/c_2)^{M+N}$. Then almost surely, for $n \in \mathbb{N}$,*
$$\mathcal{L}((L_k)_{k=1,\ldots,n} | \mathcal{F}) \leq C_1^n \mathcal{L}((r^{(4-d)k} X_k)_{k=1,\ldots,n}).$$

PROOF. As we know already the independence of $(L_k)_{k \in \mathbb{N}}$ given $\mathcal{F}$, it remains to show that
$$\mathcal{L}(L_k | \mathcal{F}) \leq C_1 \mathcal{L}(r^{(4-d)k} X_k) \qquad \text{for all } k \in \mathbb{N}.$$
By Brownian scaling the law of $r^{-(4-d)k} L_k$ given $\mathcal{F}$ is the law of $L_0$ with respect to an $(M+N)$-tuple of independent Brownian motions each started at a fixed point on $\partial B(0, 2r)$ and conditioned to exit $B(0,2)$ in a fixed point. Hence the result is a direct consequence of Lemma 2.4. $\square$

PROOF OF THEOREM 1.1, UPPER BOUND. We are now ready to prove the upper bound in Theorem 1.1.

Let $m \in \mathbb{N}$ be large enough such that, with $\theta := r^{(4-d)m}$,
$$\mathbb{P}\{X_0 < \theta\} \leq 4^{\xi_d(M,N)} r^{\xi_d(M,N)-\varepsilon}.$$
This choice is possible by the definition of $\xi_d(M,N)$ in (1.4) and by Lemma 2.5. We first look at the sequence $\delta_n = r^{(4-d)n}$, $n \in \mathbb{N}$. By Lemma 2.6 we have
$$\mathbb{P}\{\ell(B(0,1)) < \delta_n\}$$
$$\leq \mathbb{P}\left\{\sum_{k=1}^{n-m} L_k < \delta_n\right\} \leq C_1^{n-m} \mathbb{P}\left\{\sum_{k=1}^{n-m} r^{(4-d)k} X_k < \delta_n\right\}$$
$$\leq C_1^{n-m} \mathbb{P}\left\{\sum_{k=1}^{n-m} X_k < r^{(4-d)(m-n)} \delta_n\right\} = C_1^{n-m} \mathbb{P}\left\{\sum_{k=1}^{n-m} X_k < \theta\right\}$$
$$\leq (C_1 \mathbb{P}\{X_0 < \theta\})^{n-m} \leq (4^{\xi_d(M,N)} C_1)^n \theta^{-\xi_d(M,N)} \delta_n^{(\xi_d(M,N)-\varepsilon)/(4-d)}$$
$$\leq \theta^{-\xi_d(M,N)} \delta_n^{\xi_d(M,N)/(4-d)-2\varepsilon},$$

where we used (2.7) in the last inequality. Hence

$$\limsup_{n\to\infty} \frac{-\log \mathbb{P}\{\ell(B(0,1)) < \delta_n\}}{\log \delta_n} \leq -\frac{\xi_d(M,N)}{4-d} + 2\varepsilon.$$

By monotonicity, and using that $\log \delta_n / \log \delta_{n+1} \to 1$, we get the statement for arbitrary sequences $\delta \downarrow 0$. Finally, the upper bound in the assertion follows as $\varepsilon > 0$ was arbitrary. □

2.3. *Proof of Theorem* 1.1, *lower bound.* For the proof of the lower bound, in principle, we have to present one particular strategy to attain a small amount of intersection local time and then prove that this strategy is sufficiently likely.

As pointed out before, the strategy is to put no restrictions on the motions until they leave a small ball of radius $\delta$ for the first time, but demand that they do not intersect afterwards. Note, however, that paths may return to $B(0,\delta)$ and contribute to the intersection local time there by intersecting an initial piece of the path. This means that the actual decoupling at the boundary $\partial B(0,\delta)$ is rather involved. We circumvent these difficulties by following a slightly different route and give a strong upper bound for $\mathbb{P}\{\ell(B(0,\delta^{1+\varepsilon})) > \delta^{4-d}\}$ as $\delta \downarrow 0$ as well as a lower bound for $\mathbb{P}\{\ell(B(0,1) \setminus B(0,\delta^{1+\varepsilon})) = 0\}$.

Recall that our Brownian motions are stopped upon leaving $B(0,R)$ where $R \in (1,\infty)$ if $d=2$ and $R \in (1,\infty]$ if $d=3$.

LEMMA 2.7. *Let $\varepsilon > 0$. For all $\delta > 0$ sufficiently small*

$$(2.8) \qquad \mathbb{P}\{\ell(B(0,\delta^{1+\varepsilon})) > \delta^{4-d}\} \leq \exp(-\delta^{-\varepsilon/4}).$$

PROOF. By [13], Theorem 1.1, there exists a constant $\theta = \theta(R) \in (0,\infty)$ such that

$$(2.9) \qquad \lim_{a\to\infty} a^{-1/2} \log \mathbb{P}\{\ell(B(0,1)) > a\} = -\theta.$$

The probabilities on the left-hand side are increasing in $R$, the constant $\theta$ depends on $R$, but in the case $d=3$, we have $\theta(R) \downarrow \theta(\infty) > 0$ as $R \uparrow \infty$. Hence, for $d=3$ we can restrict attention to the case $R=\infty$ and get, by Brownian scaling, for sufficiently small $\delta > 0$,

$$(2.10) \quad \mathbb{P}\{\ell(B(0,\delta^{1+\varepsilon})) > \delta\} = \mathbb{P}\{\ell(B(0,1)) > \delta^{-\varepsilon}\} \leq \exp(-\delta^{-\varepsilon/4}).$$

For the case $d=2$ we have to spend a little more work, as Brownian scaling does not apply directly. We have to consider our Brownian motions $W^i$ stopped upon leaving $B(0,R)$ for different values of $R$ now and write $W_R^i$, $\ell_R$ and so on for the corresponding random objects. For $\delta > 0$ let $R(\delta) := \delta^{-(1+\varepsilon)} R$. Now Brownian scaling yields

$$\mathbb{P}\{\ell(B(0,\delta^{1+\varepsilon})) > \delta^2\} = \mathbb{P}\{\ell_R(B(0,\delta^{1+\varepsilon})) > \delta^2\} = \mathbb{P}\{\ell_{R(\delta)}(B(0,1)) > \delta^{-2\varepsilon}\}.$$

For $i = 1, \ldots, M+N$ let $\tau^i(1) := \tau^i_{R,1}$, $\tau^i(2) := \tau^i_{R,1,R,1}$ [recall the notation from (2.2)] and so on. Define

$$X^i := \min\{n \in \mathbb{N} : \tau^i(n) > \tau^i_{R(\delta)}\};$$

this is (one plus) the number of downcrossings of the annulus $B(0,R) \setminus B(0,1)$ by the stopped Brownian motion $W^i_{R(\delta)}$.

Before we continue the main argument we establish some auxiliary inequalities for the $X^i$. The distribution of $X^i$ is geometric with failure parameter

$$\begin{aligned} p(\delta) &:= \mathbb{P}\{\tau^i(1) < \tau^i_{R(\delta)}\} \\ &= \mathbb{P}_x\{\tau_1 < \tau_{R(\delta)}\} \qquad [\text{for } x \in \partial B(0,R)] \\ &= 1 - \frac{\log(R)}{\log(R(\delta))} = 1 - \frac{\log(R)}{\log(R) - (1+\varepsilon)\log(\delta)}. \end{aligned}$$

Hence

$$\mathbb{E}[X^i] = \frac{1}{1-p(\delta)} = \frac{\log(R) - (1+\varepsilon)\log(\delta)}{\log(R)},$$

and for $K \in \mathbb{N}$,

$$\mathbb{E}[X^1 \mathbb{1}_{\{X^1 \geq K\}}] = \left(\frac{K}{p(\delta)} + \frac{1}{1-p(\delta)}\right) p(\delta)^K.$$

In particular, for $K = K(\delta) = \delta^{-\varepsilon/2}$ and $\delta > 0$ sufficiently small

(2.11) $(M+N)MN\mathbb{E}[X^1]^{M+N-1}\mathbb{E}[X^1 \mathbb{1}_{\{X^1 \geq K(\delta)\}}] < \tfrac{1}{2}\exp(-\delta^{-\varepsilon/4}).$

On the other hand, by (2.9) for $\delta > 0$ small enough

(2.12) $MNK(\delta)^2 \mathbb{P}\left\{\ell_R(B(0,1)) > \dfrac{\delta^{-2\varepsilon}}{MN K(\delta)^2}\right\} \leq \dfrac{1}{2}\exp(-\delta^{-\varepsilon/4}).$

Now we come back to the main argument. A simple coupling argument shows that the contribution of $W^i$ to $\ell$ between $\tau^i(n)$ and $\tau^i(n+1)$ is stochastically no larger than the contribution between time 0 and $\tau^i(1)$, that is, between 0 and $\tau^i_R$. For $k = (k^1, \ldots, k^{M+N}) \in \mathbb{N}^{M+N}$ abbreviate

$$\sigma(k) := (k^1 + \cdots + k^M)(k^{M+1} + \cdots + k^{M+N}) \leq MN \prod_{i=1}^{M+N} k^i.$$

We get

$$\mathbb{P}\{\ell_{R(\delta)}(B(0,1)) > \delta^{-2\varepsilon}\}$$
$$\leq \sum_{k^1=1}^{\infty} \cdots \sum_{k^{M+N}=1}^{\infty} \sigma(k)\mathbb{P}\Big\{\ell_R(B(0,1)) > \frac{\delta^{-2\varepsilon}}{\sigma(k)};$$
$$X^i = k^i, \text{ for all } i=1,\ldots,M+N\Big\}.$$

Decomposing the sum into the contribution coming from $k^1,\ldots,k^{M+N}$ all smaller than $K(\delta)$ on the one hand, and the contribution coming from all $k^1,\ldots,k^{M+N}$ with some $k^i > K(\delta)$ on the other hand, we can bound the right-hand side by

$$MNK(\delta)^2 \mathbb{P}\Big\{\ell_R(B(0,1)) > \frac{\delta^{-2\varepsilon}}{MNK(\delta)^2}\Big\}$$
$$+ (M+N)MN\mathbb{E}[X^1]^{M+N-1}\mathbb{E}[X^1 \mathbb{1}_{\{X^1 \geq K(\delta)\}}]$$
$$\leq \frac{1}{2}\exp(-\delta^{-\varepsilon/4}) + \frac{1}{2}\exp(-\delta^{-\varepsilon/4}),$$

where we used (2.11) and (2.12) in the last step. This was the claim. □

The second ingredient for the proof of the lower bound is an estimate on the probability that a certain annulus has zero intersection local time.

LEMMA 2.8. *Let $R \in (1,\infty)$ if $d=2$, and $R \in (1,\infty]$ if $d=3$. Then*

$$\liminf_{\delta \downarrow 0} \frac{\log \mathbb{P}\{\ell(B(0,1) \setminus B(0,\delta)) = 0\}}{-\log \delta} \geq -\xi_d(M,N).$$

PROOF. Denote, for $R > 0$,

$$\mathcal{W}^1(R) = \bigcup_{i=1}^{M} W^i([0,\tau_R^i]) \quad \text{and} \quad \mathcal{W}^2(R) = \bigcup_{i=M+1}^{M+N} W^i([0,\tau_R^i]).$$

Further for $r > 0$ let $\mathbb{P}_r$ denote the probability measure under which $W^i(0)$, $i=1,\ldots,M+N$, are independent and uniformly distributed on $\partial B(0,r)$. In the case $R < \infty$ we have by (1.4) and Brownian scaling, as $\delta \downarrow 0$ (here we do not need Lemma 2.5 which would yield equality in the first step)

$$\mathbb{P}\{\ell(B(0,1) \setminus B(0,\delta)) = 0\} \geq \mathbb{P}_\delta\{\mathcal{W}^1(R) \cap \mathcal{W}^2(R) = \varnothing\}$$
(2.13)
$$\geq \Big(\frac{\delta}{R}\Big)^{\xi_d(M,N)+o(1)}.$$

In the case $d = 3$ and $R = \infty$ this estimate is apparently not good enough. However, due to transience, we can postulate that our Brownian motions do not return to $B(0,1)$ once they have left $B(0,2)$ and then apply (2.13) with $R = 2$. Indeed, consider the events [recall (2.2)]

$$A^i = \{\tau_{2,1}^i = \infty\} \quad \text{and} \quad A = \bigcap_{i=1}^{M+N} A^i.$$

Note that $\mathbb{P}_\delta(A^i) = \frac{1}{2}$ for all $i$ and $\delta \in [0,1]$, hence $\mathbb{P}_\delta(A) = 2^{-(M+N)}$. By the strong Markov property applied to $\tau_2^i$, the family $(\{\mathcal{W}^1(2) \cap \mathcal{W}^2(2) = \varnothing\}, A^1, \ldots, A^{M+N})$ is independent under $\mathbb{P}_\delta$ for all $\delta \in (0,1)$. Thus

$$\begin{aligned}
\mathbb{P}\{\ell(B(0,1) \setminus B(0,\delta)) = 0\} \\
\geq \mathbb{P}_\delta(\{\mathcal{W}^1(2) \cap \mathcal{W}^2(2) = \varnothing\} \cap A) \\
\geq 2^{-(M+N+\xi_d(M,N))} \delta^{\xi_d(M,N)+o(1)} \qquad \text{as } \delta \downarrow 0. \quad \square
\end{aligned}$$

PROOF OF THEOREM 1.1, LOWER BOUND. We can now assemble the pieces. Fix $\varepsilon > 0$. We make the simple observation that for $\delta > 0$

$$\begin{aligned}
\mathbb{P}\{\ell(B(0,1)) < \delta^{4-d}\} \\
\geq \mathbb{P}\{\ell(B(0,1) \setminus B(0,\delta^{1+\varepsilon})) = 0\} - \mathbb{P}\{\ell(B(0,\delta^{1+\varepsilon})) \geq \delta^{4-d}\}.
\end{aligned}$$

By Lemmas 2.7 and 2.8 the second term on the right-hand side is of smaller order than the first term, which is of order $\geq \delta^{\xi_d(M,N)(1+\varepsilon)+o(1)}$ as $\delta \downarrow 0$. This yields

$$\liminf_{\delta \downarrow 0} \frac{\log \mathbb{P}\{\ell(B(0,1)) < \delta\}}{-\log \delta} \geq -\frac{\xi_d(M,N)}{4-d}(1+\varepsilon).$$

As $\varepsilon > 0$ was arbitrary, the lower bound of Theorem 1.1 is established. This completes the proof of the theorem. $\square$

2.4. *Reversing paths*: *local versions of the tail asymptotics.* The aim of this section is to prove the two results, Lemma 2.9 and Lemma 2.11, which reformulate the tail asymptotics established in the previous sections in a form suitable for use in the proof of Theorem 1.3. The two proofs are largely analogous and make the time-reversal of paths (mentioned in the overview, Section 1.4) precise. We start with the result needed for the lower bound, which strictly speaking is a reformulation of the definition of the intersection exponents.

Recall from (2.2) that for $r, s > 0$, $\tau_r^i$ is the first hitting time of $\partial B(0,r)$ and $\tau_{r,s}^i$ is the first hitting time of $\partial B(0,s)$ after $\tau_r^i$, for the Brownian motion $W^i$.

LEMMA 2.9. *Fix $b > 1 > c$ and $r > 0$. Suppose that $W^1, W^2$ are independent Brownian paths started uniformly on the sphere $\partial B(0,r) \subset \mathbb{R}^d$, for $d = 2, 3$. Then*

$$\lim_{r \downarrow 0} \frac{1}{\log(1/r)} \log \mathbb{P}\{W^1([0, \tau^1_{r^b,r}])$$

(2.14)
$$\cap W^2([0, \tau^2_{r^b,r}]) = \varnothing \mid \tau^1_{r^b} < \tau^1_{r^c}, \tau^2_{r^b} < \tau^2_{r^c}\}$$

$$= \xi_d(2,2)(1-b).$$

PROOF. We look at the *upper bound* and define random times

$$\tau^i_* = \sup\{t < \tau^i_{r^b} : |W^i(t)| = r\} \qquad \text{for } i = 1, 2.$$

The paths $e^i : [0, \tau^i_{r^b} - \tau^i_*] \to \mathbb{R}^d$, $e^i(t) = W^i(t + \tau^i_*)$, are *Brownian excursions* from $\partial B(0, r)$ to $\partial B(0, r^b)$, and hence the time-reversed paths

$$e^i_* : [0, \tau^i_{r^b} - \tau^i_*] \to \mathbb{R}^d, \qquad e^i_*(t) = e^i(\tau^i_{r^b} - \tau^i_* - t)$$

are Brownian excursions from $\partial B(0, r^b)$ to $\partial B(0, r)$. Now fix $b > \beta > 1$ and define $\sigma^i = \inf\{t > 0 : |e^i_*(t)| = r^\beta\}$. As the transition semigroup of a Brownian excursion in $B(0,r) \setminus B(0, r^b)$ is the same as for Brownian motion killed upon leaving $B(0,r) \setminus B(0, r^b)$, the processes

$$\overline{W}^i : [0, \tau^i_{r^b} - \tau^i_* - \sigma^i] \to \mathbb{R}^d, \qquad \overline{W}^i(t) = e^i_*(\sigma^i + t),$$

are independent Brownian motions, started in a uniformly chosen point on $\partial B(0, r^\beta)$ killed upon leaving $B(0,r) \setminus B(0, r^b)$ and conditioned to hit $\partial B(0, r)$ before $\partial B(0, r^b)$. Denoting the first hitting times of $\partial B(0, s)$ by the motion $\overline{W}^i$ by $\overline{\tau}^i_s$, we get

$$\mathbb{P}\{W^1[0, \tau^1_{r^b,r}] \cap W^2[0, \tau^2_{r^b,r}] = \varnothing \mid \tau^1_{r^b} < \tau^1_{r^c}, \tau^2_{r^b} < \tau^2_{r^c}\}$$

$$\leq \mathbb{P}\{(W^1[\tau^1_*, \tau^1_{r^b}] \cup W^1[\tau^1_{r^b}, \tau^1_{r^b,r}])$$

$$\cap (W^2[\tau^2_*, \tau^2_{r^b}] \cup W^2[\tau^2_{r^b}, \tau^2_{r^b,r}]) = \varnothing \mid \tau^1_{r^b} < \tau^1_{r^c}, \tau^2_{r^b} < \tau^2_{r^c}\}$$

$$\leq \mathbb{P}\{(\overline{W}^1[0, \overline{\tau}^1_r] \cup W^1[\tau^1_{r^b}, \tau^1_{r^b,r}])$$

$$\cap (\overline{W}^2[0, \overline{\tau}^2_r] \cup W^2[\tau^2_{r^b}, \tau^2_{r^b,r}]) = \varnothing \mid \tau^i_{r^b} < \tau^i_{r^c}, \overline{\tau}^i_r < \overline{\tau}^i_{r^b} \text{ for } i = 1,2\}$$

$$\leq \mathbb{P}\{(\overline{W}^1[0, \overline{\tau}^1_r] \cup \widetilde{W}^1[0, \widetilde{\tau}^1_r])$$

$$\cap (\overline{W}^2[0, \overline{\tau}^2_r] \cup \widetilde{W}^2[0, \widetilde{\tau}^2_r]) = \varnothing \mid \overline{\tau}^1_r < \overline{\tau}^1_{r^b}, \overline{\tau}^2_r < \overline{\tau}^2_{r^b}\},$$

where $\widetilde{W}^i$ is a Brownian motion which [except for the starting point on $\partial B(0, r^\beta)$] is independent of $\overline{W}^i$ and which is stopped at the time $\widetilde{\tau}^i_r$ when it first hits $B(0,r)$. By Lemma 2.1, for each Brownian motion $\overline{W}^i$ the probability of the conditioning event $\{\overline{\tau}_r^i < \overline{\tau}_{r^b}^i\}$ is equal to $(\beta-b)/(1-b)$ in $d=2$ and, in $d=3$, it is converging to 1, as $r \downarrow 0$. In any case we find $\varepsilon > 0$ such that

$$\mathbb{P}\{\overline{\tau}_r^i < \overline{\tau}_{r^b}^i\} > \varepsilon \qquad \text{for all } i=1,2 \text{ and } 0 < r < \tfrac{1}{2}.$$

We can thus continue and find

$$\mathbb{P}\{W^1[0,\tau_{r^b,r}^1] \cap W^2[0,\tau_{r^b,r}^2] = \varnothing | \tau_{r^b}^1, \tau_{r^b}^2 < \infty\}$$
$$\leq \mathbb{P}\{(\overline{W}^1[0,\overline{\tau}_r^1] \cup \widetilde{W}^1[0,\widetilde{\tau}_r^1])$$
$$\cap (\overline{W}^2[0,\overline{\tau}_r^2] \cup \widetilde{W}^2[0,\widetilde{\tau}_r^2]) = \varnothing | \overline{\tau}_r^1 < \overline{\tau}_{r^b}^1, \overline{\tau}_r^2 < \overline{\tau}_{r^b}^2\}$$
$$\leq \varepsilon^{-2} \mathbb{P}\{(\overline{W}^1[0,\overline{\tau}_r^1] \cup \widetilde{W}^1[0,\widetilde{\tau}_r^1]) \cap (\overline{W}^2[0,\overline{\tau}_r^2] \cup \widetilde{W}^2[0,\widetilde{\tau}_r^2]) = \varnothing\}.$$

By Brownian scaling the probability of the last event is equal to the probability that for two families $\mathfrak{B}^1$, $\mathfrak{B}^2$ of two Brownian paths started on the unit sphere we observe $\mathfrak{B}^1(r^{1-\beta}) \cap \mathfrak{B}^1(r^{1-\beta}) = \varnothing$. Recalling the definition of the intersection exponent $\xi_d(2,2)$, this leads to

$$\limsup_{r \downarrow 0} \frac{1}{\log(1/r)} \log \mathbb{P}\{W^1[0,\tau_{r^b,r}^1] \cap W^2[0,\tau_{r^b,r}^2] = \varnothing | \tau_{r^b}^1 < \tau_{r^c}^1, \tau_{r^b}^2 < \tau_{r^c}^2\}$$
$$\leq \xi_d(2,2)(1-\beta).$$

Letting $\beta \uparrow b$ gives the upper bound.

For the proof of the *lower bound* we argue similarly. Let $c < \gamma < 1$. Note that

$$\mathbb{P}\{W^1[0,\tau_{r^b,r}^1] \cap W^2[0,\tau_{r^b,r}^2] = \varnothing | \tau_{r^b}^1 < \tau_{r^c}^1, \tau_{r^b}^2 < \tau_{r^c}^2\}$$
$$\geq \mathbb{P}\{W^1[0,\tau_{r^b,r}^1] \cap W^2[0,\tau_{r^b,r}^2] = \varnothing | \tau_{r^b}^1 < \tau_{r^\gamma}^1, \tau_{r^b}^2 < \tau_{r^\gamma}^2\}$$
$$\times \left[\frac{\mathbb{P}\{\tau_{r^b}^1 < \tau_{r^\gamma}^1\}}{\mathbb{P}\{\tau_{r^b}^1 < \tau_{r^c}^1\}}\right]^2.$$

By Lemma 2.1, there exists $\varepsilon > 0$ and such that the fraction is bigger than $\varepsilon$ for all $0 < r < 1/2$. We can write

$$W^i[0,\tau_{r^b,r}^i] = W^i[0,\tau_{r^b}^i] \cup W^i[\tau_{r^b}^i, \tau_{r^b,r}^i].$$

Under the new conditioning the *first path* $W^i[0,\tau_{r^b}^i]$ can be seen as part of a Brownian excursion from $\partial B(0,r^\gamma)$ to $\partial B(0,r^b)$, or, by time-reversal as part of a Brownian excursion $e_*^i$ from $\partial B(0,r^b)$ to $\partial B(0,r^\gamma)$. This excursion naturally is part of a Brownian motion $\overline{W}^i$ started in a point uniformly distributed on $\partial B(0,r^b)$ and stopped upon hitting $B(0,r^\gamma)$, say at time $\overline{\tau}_{r^\gamma}^i$. By extending the *second path* $W^1[\tau_{r^b}^1, \tau_{r^b,r}^1]$ to the right by a Brownian motion path until it hits $r^\gamma$, we see that we can replace it by a Brownian motion $\widetilde{W}^1$ started at $W^1(\tau_{r^b}^1)$ and stopped at its first hitting time of $\partial B(0, r^\gamma)$, which we denote by $\widetilde{\tau}_{r^\gamma}^1$. Hence,

$$\mathbb{P}\{W^1[0,\tau_{r^b,r}^1] \cap W^2[0,\tau_{r^b,r}^2] = \varnothing | \tau_{r^b}^1 < \tau_{r^\gamma}^1, \tau_{r^b}^2 < \tau_{r^\gamma}^2\}$$
$$\geq \mathbb{P}\{(\overline{W}^1[0,\overline{\tau}_{r^\gamma}^1] \cup \widetilde{W}^1[0,\widetilde{\tau}_{r^\gamma}^1]) \cap (\overline{W}^2[0,\overline{\tau}_{r^\gamma}^2] \cup \widetilde{W}^2[0,\widetilde{\tau}_{r^\gamma}^2]) = \varnothing\}.$$

Apart from the starting point, $\overline{W}^1$ is independent of $\widetilde{W}^1$, and $\overline{W}^2$ is independent of $\widetilde{W}^2$. As the starting points $W^1(\tau_{r^b}^1)$ and $W^2(\tau_{r^b}^2)$ are independent and uniformly distributed on $\partial B(0, r^b)$ we get

$$\liminf_{r \downarrow 0} \frac{1}{\log(1/r)} \log \mathbb{P}\{W^1[0,\tau_{r^b,r}^1] \cap W^2[0,\tau_{r^b,r}^2] = \varnothing | \tau_{r^b}^1, \tau_{r^b}^2 < \infty\}$$
$$\geq \xi_d(2,2)(\gamma - b).$$

Now letting $\gamma \uparrow 1$ gives the result. □

REMARK 2.10. In Lemma 2.9 we could fix vectors $u_1, \ldots, u_4 \in \partial B(0,1)$ with $\{u_1, u_3\} \cap \{u_2, u_4\} = \varnothing$. Replacing the starting points by fixed points $W^1(0) = ru_1$ and $W^2(0) = ru_2$ and also fixing the exit points from the ball $B(0,r)$ as $W^1(\tau_{r^b,r}^1) = u_3 r$ and $W^2(\tau_{r^b,r}^2) = u_4 r$, the result remains unchanged. Moreover, the convergence is uniform in $u_1, \ldots, u_4$ as long as the minimal distance between points of $\{u_1, u_3\}$ and $\{u_2, u_4\}$ is bounded away from zero. This can be done by a standard argument, as performed, for example, in [13], Lemma 5.7; see also Lemma 2.4.

We now formulate a version of Theorem 1.1 which represents the connection between the tail behavior of the intersection local time and the multifractal spectrum. The following lemma plays a crucial role in the proof of the upper bound of Theorem 1.3.

LEMMA 2.11. *For all $a > 4 - d$, and $\varepsilon > 0$,*

$$\lim_{r \downarrow 0} \frac{1}{-\log r} \log \mathbb{P}\{\ell(B(x,r)) < r^a | \tau^1(B(x, r^{a/(4-d)})) < 1 - \varepsilon,$$
(2.15)
$$\tau^2(B(x, r^{a/(4-d)})) < 1 - \varepsilon\}$$
$$= \xi_d(2,2)\left(1 - \frac{a}{4-d}\right).$$

PROOF. We only show the upper bound in Lemma 2.11, as the lower bound is not used in the paper, and the proof is quite similar. For notational simplicity we replace $a$ by $a(4-d)$ and let

$$p(x,r,a) := \frac{1}{-\log r}\log\mathbb{P}\{\ell(B(x,r)) < r^{(4-d)a}|\tau^1(B(x,r^a)) < 1-\varepsilon,$$
$$\tau^2(B(x,r^a)) < 1-\varepsilon\}.$$

For an upper bound, we can always replace the Brownian paths by smaller pieces, effectively making the intersection local time smaller. Hence we may start the motions $W^i$ at time $\tau^i(B(x,r))$ and stop them at time $\tau^i(B(x,r^a),\partial B(x,r))$, if this time is smaller than 1. We may assume the latter as

$$\mathbb{P}\{\tau^i(B(x,r^a),\partial B(x,r)) \geq 1|\tau^i(B(x,r^a)) < 1-\varepsilon\}$$

is decaying faster than exponentially. Let $\widetilde{W}^1,\ldots,\widetilde{W}^4$ be independent Brownian motions started in $x$ and stopped upon leaving $B(x,r)$. Let $\tilde{\ell}$ denote the intersection local time of $\widetilde{W}^1,\widetilde{W}^2$ with $\widetilde{W}^3,\widetilde{W}^4$. Let $\tilde{\varepsilon} > 0$ be small. Arguing as in the proof of Lemma 2.9 we can replace $W^1([\tau^1(B(x,r)),\tau^1(\partial B(x,r^a))])$ and $W^1([\tau^1(\partial B(x,r^a)),\tau^1(B(x,r^a),\partial B(x,r))])$ by $\widetilde{W}^1,\widetilde{W}^2$ and the analogous pieces of $W^2$ by $\widetilde{W}^3,\widetilde{W}^4$. We obtain

$$\limsup_{r\downarrow 0} p(x,r,a) \leq \limsup_{r\downarrow 0} \frac{1}{-\log r}\mathbb{P}\{\tilde{\ell}(B(x,r) \setminus B(x,r^{a-\tilde{\varepsilon}})) < r^{(4-d)a}\}.$$

We argue as in the proof of the lower bound of Theorem 1.1. Write

$$\mathbb{P}\{\tilde{\ell}(B(x,r) \setminus B(x,r^{a-\tilde{\varepsilon}})) < r^{(4-d)a}\}$$
$$\leq \mathbb{P}\{\tilde{\ell}(B(x,r)) < r^{(4-d)(a-3\tilde{\varepsilon})}\} + \mathbb{P}\{\tilde{\ell}(B(x,r^{a-\tilde{\varepsilon}})) > r^{(4-d)(a-2\tilde{\varepsilon})}\}.$$

Using Lemma 2.7 we see that the second term on the right-hand side is negligible. Hence,

$$\limsup_{r\downarrow 0} p(x,r,a) \leq \limsup_{r\downarrow 0} \frac{1}{-\log r}\mathbb{P}\{\tilde{\ell}(B(x,r)) < r^{(4-d)(a-3\tilde{\varepsilon})}\}$$
$$= \xi_d(2,2)(1-(a-3\tilde{\varepsilon})),$$

from Theorem 1.1 (and Brownian scaling). The result follows as $\tilde{\varepsilon}$ can be arbitrarily close to 0. □

**3. Multifractal spectrum.** In this section we prove Theorem 1.3 by showing the upper bound and the lower bound separately.

3.1. *Proof of Theorem* 1.3, *upper bounds.* This follows from a standard first moment method.

PROPOSITION 3.1 (Upper bound). *Almost surely, for every $a \geq 4 - d$,*

$$\dim \mathcal{T}(a) \leq (4-d) \frac{\xi_d(2,2)}{a} + (4-d) - \xi_d(2,2),$$

*where negative values of the dimension indicate that the set is empty.*

PROOF. The proof makes use of Lemma 2.11. The case $a = 4 - d$ is trivial, hence we may fix $a > 4 - d$. Now fix a cube C of unit length that has positive distance, say $\delta$, to the origin. Also fix $\varepsilon > 0$ and let $\mathcal{T}_\varepsilon(a)$ be the set of $a$-thin points $x$, such that the first visit of both motions to $x$ is before time $1 - \varepsilon$. Abbreviating $\xi := \xi_d(2,2)$, it now suffices to prove

(3.1) $$\dim(\mathcal{T}_\varepsilon(a) \cap \mathtt{C}) \leq \xi \left( \frac{4-d}{a} - 1 \right) + (4-d).$$

For $k \in \mathbb{N}$ denote by $\mathfrak{D}_k$ the set of dyadic cubes with respect to C with sidelength $2^{-k}$. For each cube $E \in \mathfrak{D}_k$ and $r > 0$ let $B(E, r)$ denote the ball of radius $r$ centered in the center of $E$.

Let $4 - d < b < a$. For $k$ large enough $\mathrm{dist}(0, B(E, 2^{-k(4-d)/a})) \geq \delta/2$ for all $E \in \mathfrak{D}_k$. Hence there exists a constant $c > 0$ such that $\mathcal{L}(W^i(\tau^i(B(E, 2^{-k(4-d)/a})))) \geq c \mathcal{U}(\partial B(E, 2^{-k(4-d)/a}))$, where $\mathcal{U}$ is the uniform distribution on the boundary of $B(E, 2^{-k(4-d)/a})$. Hence by Lemma 2.11 for $4 - d < b' < b$ and for $k \geq k_0 = k_0(b, b')$ large enough,

(3.2) $$\mathbb{P}\{\ell(B(E, \tfrac{1}{2} 2^{-k(4-d)/b})) \leq 2^{-k(4-d)} | \tau^1(B(E, 2^{-k})),$$
$$\tau^2(B(E, 2^{-k})) < 1 - \varepsilon\}$$
$$\leq 2^{k\xi((4-d)/b' - 1)},$$

for all $E \in \mathfrak{D}_k$. Let

$$\mathfrak{D}_k(b) := \{ E \in \mathfrak{D}_k : \tau^i(B(E, 2^{-k})) < 1 - \varepsilon, \text{ for } i = 1, 2,$$
$$\text{there exists } x \in S \cap E \text{ with } \ell(B(x, 2^{-k(4-d)/b})) \leq 2^{-k(4-d)} \}.$$

Then, for any $k_1 \in \mathbb{N}$, the collection $\bigcup_{k \geq k_1} \mathfrak{D}_k(b)$ is a covering of $\mathcal{T}_\varepsilon(a) \cap \mathtt{C}$.

*Case $d = 3$.* There exists a constant $C$ such that

(3.3) $$\mathbb{P}\{\tau^1(B(E, 2^{-k})) < 1 - \varepsilon, \tau^2(B(E, 2^{-k})) < 1 - \varepsilon\} \leq C 2^{-2k}.$$

By (3.3) and (3.2) for any $k \geq k_0$,

$$\mathbb{P}\{E \in \mathfrak{D}_k(b)\} \leq C 2^{k\xi((1/b')-1)} 2^{-2k},$$

for all $E \in \mathfrak{D}_k$. Thus for $\alpha \geq 0$

$$\sum_{k \geq k_1} 2^{-\alpha k} \mathbb{E}[\#\mathfrak{D}_k(b)] \leq C \sum_{k \geq k_1} 2^{-\alpha k} 2^k 2^{k\xi((1/b')-1)},$$

which is finite, if and only if $\alpha > 1 + \xi((1/b') - 1)$. This yields

$$\dim(\mathcal{T}_\varepsilon(a) \cap \mathtt{C}) \leq 1 + \xi\left(\frac{1}{b'} - 1\right).$$

As $b' \in (1, a)$ could be chosen arbitrarily close to $a$, this yields the upper bound

$$\dim(\mathcal{T}(a) \cap \mathtt{C}) \leq 1 + \xi\left(\frac{1}{a} - 1\right).$$

*Case $d = 2$.* In this case

$$\mathbb{P}\{E \in \mathfrak{D}_k(b)\} \leq \mathbb{P}\{E \in \mathfrak{D}_k(b) | \tau^1(B(E, 2^{-k})), \tau^2(B(E, 2^{-k})) < 1 - \varepsilon\}$$
$$\leq 2^{k\xi(2/b' - 1)}.$$

Thus $\mathbb{E}[\#\mathfrak{D}_k(b)] \leq 2^{2k} 2^{\xi k(2/b' - 1)}$. Continuing the argument as above yields the claim. □

3.2. *Lower bounds: the percolation technique.* In order to prove the lower bound, we fix $R > \sqrt{d}$ and work with the intersection local time $\ell$ of two Brownian motions running up to the first exit time for a large ball $B_R := B(0, R)$. We denote the set of $a$-thin points, respectively, strictly $a$-thin points, by $\mathcal{T}(a, R)$, respectively, $\mathcal{T}^s(a, R)$. The arguments following Proposition 3.5 show how to get rid of this assumption and get the bound for Brownian motions running for any finite amount of time.

To obtain lower bounds we use the method of intersection with independent random sets; see, for example, [12] for an extensive account of this. However, to realize this method new techniques are needed. Compared to the approach of [12] we are facing two additional difficulties: on the one hand the presence of long-range dependence thanks to the recurrence of Brownian motion in $d = 2$, and on the other hand the lack of a natural parametrization of $S$ by a nonrandom set. Note in particular, that $\mathcal{T}(a)$ are *not* lim sup random fractals in the sense of [12] or [2], as they are not dense in a nonrandom set. We shall show in Sections 3.3 and 3.4 how to overcome these difficulties by adapting and combining ideas of [2], which handles long-range dependence, and of [13], which deals with subfractals of random sets.

Suppose now that $\mathtt{C} \subset \mathbb{R}^d$ is a fixed compact unit cube not containing the origin. We denote by $\mathfrak{D}_n$ the collection of *compact* dyadic subcubes (relative to $\mathtt{C}$) of sidelength $2^{-n}$. We also let $\mathfrak{D} = \bigcup_{n=0}^\infty \mathfrak{D}_n$. Given $\gamma \in [0, d]$ we construct a random compact set $\Gamma[\gamma] \subset \mathtt{C}$ inductively as follows: We keep each of the $2^d$ compact cubes in $\mathfrak{D}_1$ independently with probability $p = 2^{-\gamma}$. Let $\mathcal{P}_1$ be the collection of cubes kept in this procedure and let $\Gamma_1[\gamma]$ be their union. Pass from $\mathcal{P}_n$ to $\mathcal{P}_{n+1}$ by keeping each cube of $\mathfrak{D}_{n+1}$, which is not contained in a previously rejected cube, independently with probability $p$, and again let $\Gamma_{n+1}[\gamma]$ be the union of the surviving cubes.

DEFINITION 3.2. The random set
$$\Gamma[\gamma] := \bigcap_{n=1}^{\infty} \Gamma_n[\gamma]$$
is called a percolation limit set.

The usefulness of percolation limit sets in fractal geometry comes from the following lemma (see, e.g., [23] for a proof).

LEMMA 3.3. *For every $\gamma \in [0,d]$ and every Borel set $A \subset \mathsf{C}$ the following properties hold:*

(i) *if $\dim A < \gamma$, then $\mathbb{P}\{A \cap \Gamma[\gamma] \neq \varnothing\} = 0$,*
(ii) *if $\dim A > \gamma$, then $\mathbb{P}\{A \cap \Gamma[\gamma] \neq \varnothing\} > 0$,*
(iii) *if $\dim A > \gamma$, then*
$$\mathbb{P}\{\dim(A \cap \Gamma[\gamma]) \leq \dim A - \gamma\} = 1$$
*and*
$$\mathbb{P}\{\dim(A \cap \Gamma[\gamma]) \geq \dim A - \gamma - \varepsilon\} > 0 \qquad \text{for all } \varepsilon > 0.$$

We now suppose that the random set $\Gamma[\gamma]$ and two Brownian motions $W^1$ and $W^2$, started at the origin, are realized independently on the same probability space, and we write $\mathbb{P}$ for the joint distribution of the motions and $\Gamma[\gamma]$. Observe that the first part of Lemma 3.3 gives a *lower* bound $\gamma$ for the Hausdorff dimension of the set $\mathcal{T}(a,R)$, if we can show that $\mathcal{T}(a,R) \cap \Gamma[\gamma] \neq \varnothing$ with positive probability. The following lemma shows that this approach also allows us to compare the sets of thin and strictly thin points (recall Definition 1.2). Recall that we abbreviate $\xi = \xi_d(2,2)$.

LEMMA 3.4. *If $\gamma = (4-d)\frac{\xi}{a} + (4-d) - \xi$, then*
(3.4) $$\mathbb{P}\{\mathcal{T}(a,R) \cap \Gamma[\gamma] = \mathcal{T}^s(a,R) \cap \Gamma[\gamma]\} = 1.$$

PROOF. An obvious modification of the upper bound established in Proposition 3.1 shows that $\dim(\mathcal{T}(a+\frac{1}{n},R)) < \gamma$, and, by Lemma 3.3(ii), we have that $\mathcal{T}(a+\frac{1}{n},R) \cap \Gamma[\gamma] = \varnothing$ almost surely for all $n$. Hence, almost surely,
$$\mathcal{T}^s(a,R) \cap \Gamma[\gamma] = \mathcal{T}(a,R) \cap \Gamma[\gamma] \cap \bigcap_{n=1}^{\infty} \mathcal{T}(a+1/n,R)^c$$
$$= \mathcal{T}(a,R) \cap \Gamma[\gamma]. \qquad \square$$

Hence the crucial part in establishing the lower bound on the dimension in Theorem 1.3 is the following proposition, whose proof will be given in the subsequent sections.

PROPOSITION 3.5 (Lower bound). *Let $\gamma = (4-d)\frac{\xi}{a} + (4-d) - \xi$. Then* $\mathbb{P}\{\mathcal{T}(a,R) \cap \Gamma[\gamma] \neq \varnothing\} > 0$.

PROOF OF THEOREM 1.3, LOWER BOUND. Proposition 3.5 implies the result of Theorem 1.3 by the following consideration. We use the following simple fact noted in [3], (3.2):

Let $A \subset \mathbb{R}^d$ be a fixed analytic set, and let $W$ be a Brownian motion with arbitrary starting point; then

(3.5) $$\dim(A \setminus W([0,\infty))) = \dim A \quad \textit{almost surely.}$$

Indeed, to verify (3.5), suppose that $\dim A > \alpha$. Then, by Frostman's lemma, see, for example, [8], 4.11, there exists a measure $\nu \neq 0$ on $A$ such that $\nu(B) \leq (\operatorname{diam}(B))^\alpha$ for all balls $B$. By Fubini's theorem $\mathbb{E}[\nu(W([0,\infty)))] = \int \mathbb{P}\{x \in W([0,\infty))\}\nu(dx) = 0$, and hence $\nu$ is concentrated on $A \setminus W([0,\infty))$ almost surely. Hence, $\dim(A \setminus W([0,\infty))) \geq \alpha$ almost surely, by the mass distribution principle; see, for example, [8], 4.2.

Now fix a compact unit cube $C \subset B_R$ at positive distance from the origin. By Proposition 3.5 and Lemma 3.4 we get

$$\mathbb{P}\{\mathcal{T}^s(a,R) \cap \Gamma[\gamma] \neq \varnothing\} > 0.$$

Together with Lemma 3.3(i) this implies

$$\mathbb{P}\{\dim \mathcal{T}^s(a,R) \geq \gamma\} = p(R) > 0.$$

By Brownian scaling the probability $p(R)$ does not depend on the choice of $R > 0$, hence we may write $p = p(R) > 0$. Now define events

$$D_n := \{\dim \mathcal{T}^s(a, 1/n) \geq \gamma\} \quad \text{for all } n \in \mathbb{N}.$$

By (3.5) with $A = \mathcal{T}^s(a, \frac{1}{n+1})$ and $W([0,\infty))$ replaced by $W^1([\tau^1_{1/(n+1)}, \tau^1_{1/n}]) \cup W^2([\tau^2_{1/(n+1)}, \tau^2_{1/n}])$, and the strong Markov property we have that

$$\mathbb{P}(D_{n+1} \setminus D_n) = 0 \quad \text{for all } n \in \mathbb{N}.$$

Hence

$$\mathbb{P}\left(\bigcap_{n=1}^\infty D_n\right) = \lim_{n \to \infty} \mathbb{P}(D_n) = p > 0.$$

The event $\bigcap_{n=1}^\infty D_n$ is in the germ $\sigma$-field of Brownian motion and hence, by Blumenthal's zero–one law, the probability $p$ is actually equal to 1. Now, back to the situation where the Brownian motions are running for a fixed time, we have

$$\mathbb{P}\{\dim \mathcal{T}^s(a) \geq \gamma\} \geq \mathbb{P}\left(\bigcap_{n=1}^\infty D_n\right) - \mathbb{P}\left(\{\dim \mathcal{T}^s(a) < \gamma\} \cap \bigcap_{n=1}^\infty D_n\right).$$

The first probability on the right-hand side is 1, and the second is easily seen to vanish, using again (3.5) and the strong Markov property. Together with the upper bound, already verified in Section 3.1, this completes the proof. □

3.3. *Lower bounds*: *removing long-range dependence.* In this section we give the core argument which allows us to handle long-range dependence in our problem. We shall not refer directly to the problem of thin points in order to simplify future use of this new technique. The key result is Proposition 3.8, which shows that there exists a large number of cubes $E$ of sidelength $2^{-k}$ that are in the $k$th step of the percolation, such that both Brownian motions hit $E$ but do not return to $E$ after first leaving a slightly larger cube around $E$.

In Section 3.4 we will separate the *global* random structure of the paths, which leads to the creation of these cubes, from the *local* random structure which, given the global structure, is independent for each cube $E$. This idea of separation of a local and a global level, using conditional independence at different places, is also the key to this section (see Lemma 3.10) though it is used here on a larger scale.

As in [2] the proof of Proposition 3.8 works essentially in two different scales. In the coarse scale we use a dimension argument to make sure that there exist enough cubes of a certain type of sidelength $2^{-K}$ for some $K \ll k$. To that end we construct a subset $S^* \subset S$ with nicer regularity features. The set $S^*$ is nonempty with positive probability and all the statements in this section which hold with positive probability actually hold almost surely on the event $\{S^* \neq \varnothing\}$.

LEMMA 3.6 (Regularization). *There exists a compact set $S^* \subset S \cap \Gamma[\gamma]$ such that, almost surely, for every open set $U \subset \mathtt{C}$:*

(i) $U \cap S^* \neq \varnothing$ *implies* $\dim(U \cap S^*) > 0$,
(ii) $\dim(U \cap S \cap \Gamma[\gamma]) > 0$ *implies* $U \cap S^* \neq \varnothing$.

*Property* (ii) *implies, in particular, that* $\mathbb{P}\{S^* \neq \varnothing\} > 0$.

PROOF. To construct the set $S^*$ we fix a countable base $\mathfrak{B}$ of open subsets of $\mathtt{C}$. We define a compact random set

$$S^* = (S \cap \Gamma[\gamma]) \setminus \bigcup \{B \in \mathfrak{B} : \dim(B \cap S \cap \Gamma[\gamma]) = 0\}.$$

Clearly, it suffices to verify (i), (ii) for a fixed set $U \in \mathfrak{B}$. Suppose first that $U \cap S^* \neq \varnothing$; then $\dim(U \cap S \cap \Gamma[\gamma]) > 0$ and hence $\dim(U \cap S^*) > 0$, which establishes (i). If $\dim(U \cap S \cap \Gamma[\gamma]) > 0$, then $U \cap S^* \neq \varnothing$ by construction, which shows (ii). □

For $U \subset \mathbb{R}^d$ let
$$\mathfrak{D}_k(U) := \{E \in \mathfrak{D}_k : E \subset U\}.$$

Now fix a bounded open set $U$ and note that there exists a constant $c(U) \in (0, \infty)$ such that

(3.6) $\quad c(U)^{-1} 2^{dk} \leq \#\mathfrak{D}_k(U) \leq c(U) 2^{dk} \quad$ for all $k \in \mathbb{N}$.

For $k \in \mathbb{N}$ and $i = 1, 2$ consider the set of cubes

(3.7) $\begin{aligned}\mathcal{H}_k^i &:= \{E \in \mathfrak{D}_k(U) : \tau^i(E) < \tau^i(B_R^c)\}, \\ \mathcal{H}_{k,\delta}^i &:= \{E \in \mathfrak{D}_k(U) : \tau^i(B(E, 2^{-k(1-\delta)})) < \tau^i(B_R^c)\} \quad \text{for } \delta \in [0,1),\end{aligned}$

that are hit by the $i$th motion, respectively, where a certain ball around the box is hit. We also write

(3.8) $\quad \mathcal{H}_k := \mathcal{P}_k \cap \mathcal{H}_k^1 \cap \mathcal{H}_k^2 \quad \text{and} \quad \mathcal{H}_{k,\delta} := \mathcal{P}_k \cap \mathcal{H}_{k,\delta}^1 \cap \mathcal{H}_{k,\delta}^2.$

DEFINITION 3.7 (*Admissible cubes*). Fix $\varepsilon > 0$ and consider the subset of those cubes that are hit by the $i$th motion but which are not visited again after first leaving $B(E, 2^{-(1-\varepsilon)k})$:
$$\mathcal{A}_k^i := \{E \in \mathcal{H}_k^i : \tau^i(E, B(E, 2^{-(1-\varepsilon)k})^c, E) > \tau^i(B_R^c)\}.$$

Now we define
$$\mathcal{A}_k := \mathcal{P}_k \cap \mathcal{A}_k^1 \cap \mathcal{A}_k^2$$
to be the set of admissible cubes $E \in \mathfrak{D}_k(U)$.

PROPOSITION 3.8. *Fix $\varepsilon > 0$ and let $(a_k)_{k \in \mathbb{N}}$ be a sequence of nonnegative real numbers such that $k^2 a_k \to 0$ if $d = 2$, and $a_k \to 0$ if $d = 3$. Then*
$$\lim_{k \to \infty} \mathbb{P}\{\#\mathcal{A}_k \leq a_k 2^{(4-d-\gamma)k} | U \cap S^* \neq \varnothing\} = 0.$$

The remainder of this section is devoted to the proof of Proposition 3.8. By the preceding lemma $U \cap S^* \neq \varnothing$ implies $\dim(U \cap S^*) > 0$ almost surely. Hence, it is enough to show for every $\delta \in (0, 2 - \gamma)$

(3.9) $\quad \lim_{k \to \infty} \mathbb{P}\{\#\mathcal{A}_k \leq a_k 2^{(4-d-\gamma)k} | \dim(U \cap S^*) > \delta\} = 0.$

Let $\varepsilon_0 > 0$ be arbitrary. Below we fix $K \in \mathbb{N}$ representing the coarse scale, and divide $\mathfrak{D}_K(U)$ into finitely many, say $m$, subgrids $\mathfrak{D}_K(U, 1), \ldots, \mathfrak{D}_K(U, m)$ such that

(3.10) $\quad d_\infty(V, V') \geq 4 \cdot 2^{-K} \quad$ for all $V, V' \in \mathfrak{D}_K(U, j), \ j = 1, \ldots, m,$

where we denote by $d_\infty(V,V')$ the maximum norm distance of the centers of $V$ and $V'$.

The idea is to show that there exists a constant $\tilde\varepsilon$, independent of $K$, such that $\#\{E \in \mathcal{A}_k : E \subset V\}$ is large with probability at least $\tilde\varepsilon$, for any given $V \in \mathcal{P}_K \cap \mathcal{H}_K^1 \cap \mathcal{H}_K^2$. Further one needs to show that for any $M$ one can choose $K$ so large that with high probability, there are at least $M$ such cubes $V$ in at least one subgrid $\mathfrak{D}_K(U,j)$. Finally, using some kind of independence between the blocks $V$ in $\mathfrak{D}_K(U,j)$ we infer that the left-hand side of (3.9) is at most $(1-\tilde\varepsilon)^M$. As $\tilde\varepsilon$ is independent of $M$, we can let $M$ tend to infinity to infer the statement.

Define
$$N_k(j) := \#(\mathfrak{D}_k(U,j) \cap \mathcal{H}_{k,0}).$$

By definition of the Hausdorff dimension,
$$\{\dim(U \cap S^*) > \delta\} \subset \left\{\max_{j=1}^m N_k(j) \geq 2^{\delta k} \text{ for all but finitely many } k\right\},$$

hence we get that there exists a $K = K(\varepsilon_0, \delta)$ such that $\mathbb{P}\{A_K | U \cap S^* \neq \varnothing\} \geq 1 - \varepsilon_0$, where
$$A_K := \bigcap_{k \geq K} \bigcup_{j=1}^m \{N_k(j) \geq 2^{\delta k}\},$$

and such that $2^{\delta K} \geq M$. Fix this $K$ from now on.

The next task is to impose a localization that produces the desired independence. This proceeding was inspired by ideas of [2]. Assume that $k \in \mathbb{N}$ is larger than $2^K$. For $V \in \mathfrak{D}_K(U)$ let
$$\mathcal{H}_k^i(V) := \{E \in \mathfrak{D}_k(V) : \tau^i(E) < \tau^i(V, B(V, 2^{-K+1})^c)\},$$
(3.11) $\quad \mathcal{A}_k^i(V) := \{E \in \mathcal{H}_k^i(V) : \tau^i(E, B(E, 2^{-(1-\varepsilon)k})^c, E) > \tau^i(B_R^c)\},$
$$\mathcal{A}_{k,\mathrm{loc}}^i(V) := \{E \in \mathcal{H}_k^i(V) : \tau^i(E, B(E, 2^{-(1-\varepsilon)k})^c, E) > \tau^i(B(E, 1/k)^c)\}.$$

Finally let
$$\mathcal{H}_k(V) := \mathcal{P}_k \cap \mathcal{H}_k^1(V) \cap \mathcal{H}_k^2(V),$$
(3.12) $\quad\quad \mathcal{A}_k(V) := \mathcal{P}_k \cap \mathcal{A}_k^1(V) \cap \mathcal{A}_k^2(V),$
$$\mathcal{A}_{k,\mathrm{loc}}(V) := \mathcal{P}_k \cap \mathcal{A}_{k,\mathrm{loc}}^1(V) \cap \mathcal{A}_{k,\mathrm{loc}}^2(V).$$

Clearly $\#\mathcal{A}_{k,\mathrm{loc}}(V) \geq \#\mathcal{A}_k(V)$. Also note that
$$\#\mathcal{A}_k \geq \#\mathcal{A}_k(V) \text{ for all } V \in \mathfrak{D}_K(U).$$

Note that the information about the value of $\#\mathcal{A}_{k,\mathrm{loc}}(V)$ is contained in $\mathcal{P}_k$ and in the Brownian motion paths between $\tau^i(V)$ and $\tau^i(V, B(V, 2^{-K+1})^c)$.

By construction, see (3.10), the intervals $(\tau^i(V), \tau^i(V, B(V, 2^{-K+1})^c))$ are disjoint for different $V \in \mathfrak{D}_K(U, j)$ and fixed $j$. This will later provide the necessary independence.

Let $(b_k)_{k \in \mathbb{N}}$ be a sequence of nonnegative real numbers such that $k^2 b_k \longrightarrow 0$ and $(k^3/\log(k))b_k \longrightarrow \infty$ if $d = 2$, or $b_k \longrightarrow 0$ and $(2^k/k)b_k \longrightarrow \infty$ if $d = 3$. Then we have

$$\mathbb{P}\{\#\mathcal{A}_k \leq a_k 2^{(4-d-\gamma)k} | U \cap S^* \neq \varnothing\}$$
$$\leq \mathbb{P}\{\#\mathcal{A}_k(V) \leq a_k 2^{(4-d-\gamma)k} \text{ for all } V \in \mathfrak{D}_K(U) | U \cap S^* \neq \varnothing\}$$
$$\leq \mathbb{P}\{U \cap S^* \neq \varnothing\}^{-1}$$
$$\times \sum_{V \in \mathfrak{D}_K(U)} \mathbb{P}[\{\#\mathcal{A}_{k,\text{loc}}(V) - \#\mathcal{A}_k(V) > b_k 2^{(4-d-\gamma)k}\} \cap A_K]$$
(3.13)
$$+ \mathbb{P}\{U \cap S^* \neq \varnothing\}^{-1}$$
$$\times \mathbb{P}[\{\#\mathcal{A}_{k,\text{loc}}(V) \leq (a_k + b_k) 2^{(4-d-\gamma)k}$$
$$\text{for all } V \in \mathfrak{D}_K(U)\} \cap A_K]$$
$$+ \mathbb{P}\{A_K^c | U \cap S^* \neq \varnothing\}$$
$$=: \mathbb{P}\{U \cap S^* \neq \varnothing\}^{-1}(I_k^1 + I_k^2) + \mathbb{P}\{A_K^c | U \cap S^* \neq \varnothing\}.$$

As $\mathbb{P}\{A_K^c | U \cap S^* \neq \varnothing\} \leq \varepsilon_0$, it suffices to show that $\limsup I_k^1 = 0$ and $\limsup I_k^2 \leq (1 - \tilde{\varepsilon})^M$ for some $\tilde{\varepsilon}$ independent of $M$.

*Estimate of $I_k^1$.* We do this estimate by first computing the first moment of $\#\mathcal{A}_{k,\text{loc}}(V)$. As we need it again later, we formulate the result as a lemma.

LEMMA 3.9. *There are constants $c_1, c_2 \in (0, \infty)$ depending only on $U$ such that for $x^1, x^2 \in \partial B(V, 2^{-K})$ and for $k \geq 2K$,*

$$c_1 \left(\frac{\varepsilon}{k}\right)^2 2^{(2-\gamma)(k-K)} \leq \mathbb{E}_{x^1,x^2}[\#\mathcal{A}_{k,\text{loc}}(V)] \leq c_2 \left(\frac{\varepsilon}{k}\right)^2 2^{(2-\gamma)(k-K)} \quad \text{if } d = 2,$$

*and*

$$c_1 2^{(1-\gamma)(k-K)} \leq \mathbb{E}_{x^1,x^2}[\#\mathcal{A}_{k,\text{loc}}(V)] \leq c_2 2^{(1-\gamma)(k-K)} \quad \text{if } d = 3,$$

*where we used the abbreviation*

$$\mathbb{E}_{x^1,x^2}[\#\mathcal{A}_{k,\text{loc}}(V)]$$
$$:= \mathbb{E}[\#\mathcal{A}_{k,\text{loc}}(V) | V \in \mathcal{H}_{K,0}, W^i(\tau^i(B(V, 2^{-K}))) = x^i \text{ for } i = 1, 2].$$

PROOF. We formulate the proof only for the upper bound in the case $d = 2$. The other three cases are quite similar. By Corollary 2.2 there exists an absolute constant $\tilde{c}$ such that for all $y^i \in \partial B(E, 2^{-(1-\varepsilon)k})$ and all $k$,

$$\mathbb{P}_{y^i}\{\tau(E) > \tau(B(E, 1/k)^c)\} \leq \tilde{c} \frac{\varepsilon k}{k - \log(k)/\log(2)} \leq 2\varepsilon\tilde{c},$$

where $\mathbb{P}_{y^i}$ and $\tau$ refer to a Brownian motion started in $y^i$. We write $\mathbb{P}_{x^1,x^2}$ for the probability measure corresponding to $\mathbb{E}_{x^1,x^2}$. Again by Corollary 2.2 and the strong Markov property of Brownian motion applied at $\tau^i(B(V, 2^{-K}))$ we get

$$\mathbb{P}_{x^1,x^2}\{\tau^i(E) < \tau^i(V, B(V, 2^{-K+1})^c)\} \leq \frac{\tilde{c}}{k - K + 1} \leq \frac{2\tilde{c}}{k}.$$

Thus, by independence and the strong Markov property of Brownian motion applied at $\tau^i(E, B(E, 2^{(1-\varepsilon)k})^c)$, we get

$$\mathbb{E}_{x^1,x^2}[\#\mathcal{A}_{k,\mathrm{loc}}(V)]$$
$$= \mathbb{E} \sum_{E \in \mathfrak{D}_k(V) \cap \mathcal{P}_k} \mathbb{P}_{x^1,x^2}\{\tau^i(E) < \tau^i(V, B(V, 2^{-K+1})^c),$$
$$\tau^i(E, B(E, 2^{-(1-\varepsilon)k})^c, E) > \tau^i(B(E, 1/k)^c)$$
$$\text{for all } i = 1, 2\}$$
$$\leq 4\tilde{c}^2\varepsilon^2 \mathbb{E} \sum_{E \in \mathfrak{D}_k(V) \cap \mathcal{P}_k} \mathbb{P}_{x^1,x^2}\{\tau^i(E) < \tau^i(V, B(V, 2^{-K+1})^c) \text{ for all } i = 1, 2\}$$
$$\leq 16\tilde{c}^4\varepsilon^2 k^{-2} \mathbb{E}[\#\{E \in \mathfrak{D}_k(V) \cap \mathcal{P}_k\} | V \in \mathcal{P}_K]$$
$$= 16\tilde{c}^4\varepsilon^2 k^{-2} 2^{(2-\gamma)(k-K)}.$$

This yields the upper bound for $d = 2$ with $c_2 = 16\tilde{c}^4$. □

The next step is to compute $\mathbb{E}[\#\mathcal{A}_{k,\mathrm{loc}}(V) - \#\mathcal{A}_k(V)]$. By definition

$$\#\mathcal{A}_{k,\mathrm{loc}}(V) - \#\mathcal{A}_k(V)$$
$$= \#(\mathcal{A}_{k,\mathrm{loc}}(V) \setminus \mathcal{A}_k(V))$$
$$= \#(\mathcal{P}_k \cap \mathcal{A}^1_{k,\mathrm{loc}}(V) \cap \mathcal{A}^2_{k,\mathrm{loc}}(V) \cap ((\mathcal{A}^1_k(V))^c \cup (\mathcal{A}^2_k(V))^c)).$$

For given $E \in \mathfrak{D}_k(V)$, by independence of the Brownian motions and the percolation,

$$\mathbb{P}\{E \in \mathcal{A}_{k,\mathrm{loc}}(V)\}$$
(3.14)
$$= 2^{-\gamma k} \prod_{i=1}^{2} \mathbb{P}\{\tau^i(E) < \tau^i_R,$$
$$\tau^i(E, B(E, 1/k)^c) < \tau^i(E, B(E, 2^{-(1-\varepsilon)k})^c, E)\}.$$

Note that $\tau^i(E, B(E, 1/k)^c) = \tau^i(E, B(E, 2^{-(1-\varepsilon)k})^c, B(E, 1/k)^c)$. Hence by Corollary 2.2 and the strong Markov property applied to $\tau^i(E, B(E, 2^{-(1-\varepsilon)k})^c)$ we can bound (for $k \geq k_0$ depending only on $R$) the right-hand side of (3.14) by

$$2^{-\gamma k} \tilde{c}^4 \left(\frac{\log(R)}{k\log(2) + \log(R)}\right)^2 \left(\frac{\varepsilon k}{k - \log(k)/\log(2)}\right)^2$$
$$\leq 2\tilde{c}^4 \varepsilon^2 2^{-\gamma k} k^{-2} \qquad \text{if } d = 2,$$

and similarly by

$$2\tilde{c}^4 2^{-\gamma k} 2^{-2k} \qquad \text{if } d = 3.$$

Using again the strong Markov property with $\tau^i(E, B(E, 1/k)^c)$ we obtain for all $E \in \mathfrak{D}(V)$, $i = 1, 2$ and $k \geq k_0$ depending only on $R$ and $\varepsilon$,

$$\mathbb{P}\{E \in (\mathcal{A}_k^i(V))^c | E \in \mathcal{A}_{k,\text{loc}}(V)\} = \mathbb{P}\{E \in (\mathcal{A}_k^i(V))^c | E \in \mathcal{A}_{k,\text{loc}}^i(V)\}$$
$$\leq \sup_{x \in \partial B(E, 1/k)} \mathbb{P}_x\{\tau(E) < \tau_R\}$$
$$\leq \begin{cases} \tilde{c}\dfrac{\log k + \log R}{k\log 2 + \log R} \leq 2\tilde{c}k^{-1}\log(k), & \text{if } d = 2, \\ \tilde{c}\dfrac{Rk - 1}{R2^k - 1} \leq \tilde{c}2^{-k}k, & \text{if } d = 3. \end{cases}$$

Altogether we have for $k \geq k_0$ and $E \in \mathfrak{D}_k(V)$

$$\mathbb{P}\{E \in \mathcal{A}_{k,\text{loc}}(V) \setminus \mathcal{A}_k(V)\} \leq \begin{cases} 4\tilde{c}^5 \varepsilon^2 2^{-\gamma k} k^{-3} \log(k), & \text{if } d = 2, \\ 4\tilde{c}^5 2^{-\gamma k} 2^{-3k} k, & \text{if } d = 3. \end{cases}$$

As $\#\mathfrak{D}_k(V) = 2^{d(k-K)}$ we conclude that there is a constant $C$ such that for $k \geq k_0$

$$(3.15) \quad \mathbb{E}[\#\mathcal{A}_{k,\text{loc}}(V) - \#\mathcal{A}_k(V)] \leq \begin{cases} C2^{(2-\gamma)k} k^{-3} \log(k), & \text{if } d = 2, \\ C2^{(1-\gamma)k} 2^{-k} k, & \text{if } d = 3. \end{cases}$$

Hence by the definition of $I_k^1$ and by the assumptions made on $(b_k)$ we get by Markov's inequality,

$$I_k^1 \leq \#\mathfrak{D}_K(U)(b_k^{-1} 2^{-(4-d-\gamma)k}) \sup_{V \in \mathfrak{D}_K(U)} \mathbb{E}[\#\mathcal{A}_{k,\text{loc}}(V) - \#\mathcal{A}_k(V)]$$
$$\leq C \times \begin{cases} \dfrac{\log(k)}{k^3 b_k}, & \text{if } d = 2 \\ k2^{-k}/b_k, & \text{if } d = 3 \end{cases} \longrightarrow 0 \qquad \text{as } k \to \infty.$$

*Estimate of $I_k^2$.* The estimate of the second term in (3.13) is much more delicate. Clearly

$$I_k^2 \leq \max\{\mathbb{P}\{\#\mathcal{A}_{k,\mathrm{loc}}(V) \leq (a_k + b_k)2^{(4-d-\gamma)k}$$
$$\text{for all } V \in \mathfrak{D}_K(U,j),\, N_K(j) \geq M\},\, j = 1, \ldots, m\}.$$

Now fix $j$. The strategy is to introduce a $\sigma$-field $\mathcal{F}_K(j)$ such that the events $\{\#\mathcal{A}_{k,\mathrm{loc}}(V) \leq (a_k + b_k)2^{(4-d-\gamma)k},\, V \in \mathfrak{D}_K(U,j)\}$ become independent given $\mathcal{F}_K(j)$ and such that for those $V$ that contribute to $N_K(j)$, and for $k$ large enough

$$(3.16) \quad \mathbb{P}\{\#\mathcal{A}_{k,\mathrm{loc}}(V) > (a_k + b_k)2^{(4-d-\gamma)k}|\mathcal{F}_K(j)\} \geq \tilde{\varepsilon}\mathbb{1}_{\{V \in \mathfrak{D}_K(U,j) \cap \mathcal{H}_{K,0}\}}$$

for some $\tilde{\varepsilon}$ that does not depend on $K$ or $M$. One can then conclude that $\limsup_{k \to \infty} I_k^2 \leq (1-\tilde{\varepsilon})^M$, which can be made arbitrarily small by letting $M \to \infty$.

We construct a decoupling $\sigma$-algebra $\mathcal{F}_K(j)$ as in the proof of Proposition 2.3. For the moment we suppress $j$ in the notation. Formally for $i = 1, 2$ we introduce a sequence of (random) sets $(V^i(n) : n = 1, \ldots, \nu^i)$ and stopping times

$$0 = \sigma^i(0) < \rho^i(1) < \sigma^i(1) < \rho^i(2) < \cdots < \sigma^i(\nu^i) < \tau_R^i < \rho^i(\nu^i + 1),$$

by

$$\rho^i(1) := \inf\{\tau^i(B(V, 2^{-K})) : V \in \mathfrak{D}_K(U,j)\},$$
$$\tau^i(B(V^i(n), 2^{-K})) = \rho^i(n) \quad \text{[this defines } V^i(n)],$$
$$\sigma^i(n) = \tau^i(B(V^i(n), 2^{-K}), B(V^i(n), 2^{-K+1})^c)$$

(3.17)
$$\text{if } \rho^i(n) < \infty,$$
$$\rho^i(n+1) = \inf\{\tau^i(B(V, 2^{-K})) :$$
$$V \in \mathfrak{D}_K(U,j) \setminus \{V^i(1), \ldots, V^i(n)\}\},$$
$$\nu^i := \max\{n : \rho^i(n) < \tau_R^i\}.$$

Now define

$$\mathcal{F}_K^i(j) := \sigma(W^i(t + \sigma^i(n)),\, t \in [0, \rho^i(n+1) - \sigma^i(n)],\, n = 0, \ldots, \nu^i)$$

and

$$\mathcal{F}_K(j) := \mathcal{F}_K^1(j) \vee \mathcal{F}_K^2(j) \vee \sigma(\mathcal{P}_K).$$

The following lemma is immediate from the construction of $\mathcal{F}_K(j)$ and the fact that (for fixed $j$) the balls $B(V, 2^{-K+1})$, $V \in \mathfrak{D}_K(U,j)$, are mutually disjoint by (3.10).

LEMMA 3.10. *The family of random variables* $(\#\mathcal{A}_{k,\mathrm{loc}}(V),\, V \in \mathfrak{D}_K(U,j))$ *is independent conditional on* $\mathcal{F}_K(j)$.

We use the notation
$$\mathbb{P}_{\mathcal{F}_K(j)} := \mathbb{P}\{\cdot\,|\mathcal{F}_K(j)\}$$

and $\mathbb{E}_{\mathcal{F}_K(j)}$ for the corresponding conditional expectation. Hence, by Lemma 3.9 there exist constants $c_1, c_2 \in (0,\infty)$ such that for $k \geq K$ almost surely

$$c_1\left(\frac{\varepsilon}{k}\right)^2 2^{(2-\gamma)k}\mathbb{1}_{\{V\in\mathcal{H}_{K,0}\}} \leq \mathbb{E}_{\mathcal{F}_K(j)}[\#\mathcal{A}_{k,\mathrm{loc}}(V)]$$
$$\leq c_2\left(\frac{\varepsilon}{k}\right)^2 2^{(2-\gamma)k}\mathbb{1}_{\{V\in\mathcal{H}_{K,0}\}} \quad \text{if } d=2,$$

and

$$c_1 2^{(1-\gamma)k}\mathbb{1}_{\{V\in\mathcal{H}_{K,0}\}} \leq \mathbb{E}_{\mathcal{F}_K(j)}[\#\mathcal{A}_{k,\mathrm{loc}}(V)] \leq c_2 2^{(1-\gamma)k}\mathbb{1}_{\{V\in\mathcal{H}_{K,0}\}} \quad \text{if } d=3.$$

By the assumption that $k^2(a_k + b_k) \to 0$ if $d=2$, and $a_k + b_k \to 0$ if $d=3$, we get that for $k$ large enough

$$\mathbb{E}_{\mathcal{F}_K(j)}[\#\mathcal{A}_{k,\mathrm{loc}}(V)] \geq 2(a_k + b_k)2^{(4-d-\gamma)k}\mathbb{1}_{\{V\in\mathcal{H}_{K,0}\}}.$$

Thus if we can replace $\#\mathcal{A}_{k,\mathrm{loc}}(V)$ in (3.16) by $\mathbb{E}_{\mathcal{F}_K(j)}[\#\mathcal{A}_{k,\mathrm{loc}}(V)]$, then we are done. To this end we have to show tightness of $(\#\mathcal{A}_{k,\mathrm{loc}}(V)/\mathbb{E}_{\mathcal{F}_K(j)}[\#\mathcal{A}_{k,\mathrm{loc}}(V)])_{k\in\mathbb{N}}$. We do so by computing second moments. Note that for $k$ large enough on $\{V \in \mathcal{H}_{K,0}\}$, using the Paley–Zygmund inequality in the second step,

$$\mathbb{P}_{\mathcal{F}_K(j)}\{\#\mathcal{A}_{k,\mathrm{loc}}(V) > (a_k + b_k)2^{(4-d-\gamma)k}\}$$
$$\geq \mathbb{P}_{\mathcal{F}_K(j)}\left\{\#\mathcal{A}_{k,\mathrm{loc}}(V) \geq \frac{1}{2}\mathbb{E}_{\mathcal{F}_K(j)}[\#\mathcal{A}_{k,\mathrm{loc}}(V)]\right\}$$
$$\geq \frac{1}{4}\frac{\mathbb{E}_{\mathcal{F}_K(j)}[\#\mathcal{A}_{k,\mathrm{loc}}(V)]^2}{\mathbb{E}_{\mathcal{F}_K(j)}[(\#\mathcal{A}_{k,\mathrm{loc}}(V))^2]}.$$

Hence the proof of Proposition 3.8 is accomplished if we can show the following lemma.

LEMMA 3.11. *There exists a constant* $c \geq 1$ *(independent of $M$) such that, for all $j = 1,\ldots,m$ and all $V \in \mathfrak{D}_K(U,j)$*

$$\mathbb{E}_{\mathcal{F}_K(j)}[(\#\mathcal{A}_{k,\mathrm{loc}}(V))^2] \leq c\mathbb{E}_{\mathcal{F}_K(j)}[\#\mathcal{A}_{k,\mathrm{loc}}(V)]^2.$$

Before we prove the lemma, we show how Proposition 3.8 can be inferred. Clearly [recall $N_K(j) = \#(\mathfrak{D}_K(U,j) \cap \mathcal{H}_{K,0})$] with $\tilde{\varepsilon} := 1/(4c)$ we get from Lemma 3.10

$$\limsup_{k\to\infty} I_k^2$$

$$\leq \max_{j=1}^{m} \limsup_{k\to\infty} \mathbb{P}\{N_K(j) \geq M \text{ and } \#\mathcal{A}_{k,\mathrm{loc}}(V) \leq (a_k + b_k)2^{(4-d-\gamma)k}$$

$$\text{for all } V \in \mathfrak{D}_K(U,j) \cap \mathcal{H}_{K,0}\}$$

$$= \max_{j=1}^{m} \limsup_{k\to\infty} \mathbb{E}\Bigg[\Bigg(\prod_{V \in \mathfrak{D}_K(U,j) \cap \mathcal{H}_{K,0}} \mathbb{P}_{\mathcal{F}_K(j)}\{\#\mathcal{A}_{k,\mathrm{loc}}(V)$$

$$\leq (a_k + b_k)2^{(4-d-\gamma)k}\}\Bigg)$$

$$\times \mathbb{1}_{\{N_K(j) \geq M\}}\Bigg]$$

$$\leq (1-\tilde{\varepsilon})^M \longrightarrow 0 \qquad \text{as } M \to \infty.$$

PROOF OF LEMMA 3.11. We do the proof explicitly only for $d = 2$ as the case $d = 3$ is quite similar. The only difference is that we have to plug in the other hitting estimates from Corollary 2.2.

By Lemma 3.9 it is enough to show that there exists a constant $C < \infty$ that is independent of $K$ and such that for $k$ large enough

$$\mathbb{E}_{\mathcal{F}_K(j)}[(\#\mathcal{A}_{k,\mathrm{loc}}(V))^2] \leq Ck^{-4}4^{(2-\gamma)(k-K)}.$$

Let $E, F \in \mathfrak{D}_k(V)$ and let $l = 2^k d_\infty(E,F)$. Recall that $d_\infty$ is the maximum distance of the centers of $E$ and $F$. Clearly

$$\mathbb{P}_{\mathcal{F}_K(j)}\{E \in \mathcal{P}_k\} = 2^{-\gamma(k-K)}\mathbb{1}_{\{V \in \mathcal{P}_K\}}.$$

In order to compute $\mathbb{P}_{\mathcal{F}_K(j)}\{E, F \in \mathcal{P}_k\}$ we define the *genealogical distance* of $E$ and $F$

$$d_{\mathrm{gen}}(E,F) := k - \sup\{s \in \{0,\ldots,k\} : E, F \in W \text{ for some } W \in \mathfrak{D}_s\}.$$

Note that $2^{d_{\mathrm{gen}}(E,F)} \geq 2^k d_\infty(E,F) = l$. Hence, on $\{V \in \mathcal{P}_K\}$,

$$(3.18) \qquad \mathbb{P}_{\mathcal{F}_K(j)}\{E, F \in \mathcal{P}_k\} = 2^{-\gamma(k-K+d_{\mathrm{gen}}(E,F))} \leq 2^{-\gamma(k-K)}l^{-\gamma}.$$

Now we come to the hitting estimates. Assume $l = 2^k d_\infty(E,F) \geq 2$. Hence, for $i = 1, 2$, by the strong Markov property and Corollary 2.2, for $k$ large enough on $\{V \in \mathcal{H}_{K,0}\}$,

$$\mathbb{P}_{\mathcal{F}_K(j)}\{\tau^i(E) < \tau^i(V, B(V, 2^{-K+1})^c), \ \tau^i(F) < \tau^i(V, B(V, 2^{-K+1})^c)\}$$

$$\leq \mathbb{P}_{\mathcal{F}_K(j)}\{\tau^i(E) < \tau^i(V, B(V, 2^{-K+1})^c)\}$$
$$\times \sup_{x \in \partial E} \mathbb{P}_x\{\tau(F) < \tau(V, B(V, 2^{-K+1})^c)\}$$
$$+ \mathbb{P}_{\mathcal{F}_K(j)}\{\tau^i(F) < \tau^i(V, B(V, 2^{-K+1})^c)\}$$
$$\times \sup_{x \in \partial F} \mathbb{P}_x\{\tau(E) < \tau(V, B(V, 2^{-K+1})^c)\}$$
$$\leq 2\tilde{c}^2 \frac{1}{k+K-1} \frac{\log(l2^{-k}) - \log(2^{-K+1})}{\log(2^{-k}) - \log(2^{-K+1})}$$
$$\leq Ck^{-1}\left(1 - \frac{\log(l)/\log(2)}{k-K+1}\right),$$

and in particular

$$\mathbb{P}_{\mathcal{F}_K(j)}\{\tau^i(E) < \tau^i(V, B(V, 2^{-K})^c)\} \leq Ck^{-1},$$

for some constant $C$ that does not depend on $M$ or $K$. Combining these estimates we get on $\{V \in \mathcal{H}_{K,0}\}$,

$$\mathbb{E}_{\mathcal{F}_K(j)}[(\#\mathcal{A}_{k,\text{loc}}(V))^2]$$
$$\leq \sum_{E \in \mathfrak{D}_k(V)} \sum_{F \in \mathfrak{D}_k(V)} \mathbb{P}_{\mathcal{F}_K(j)}\{E, F \in \mathcal{P}_k\}$$
$$\times \mathbb{P}_{\mathcal{F}_K(j)}\{\tau^i(E) < \tau^i(V, B(V, 2^{-K+1})^c),$$
$$\tau^i(F) < \tau^i(V, B(V, 2^{-K+1})^c) \text{ for all } i = 1, 2\}$$
$$\leq \sum_{E \in \mathfrak{D}_k(V)} \sum_{l=1}^{2^{k-K}} \sum_{\substack{F \in \mathfrak{D}_k(V) \\ d_\infty(E,F) = l2^{-k}}} \mathbb{P}_{\mathcal{F}_K(j)}\{E, F \in \mathcal{P}_k\} C^2 k^{-2}\left(1 - \frac{\log(l)/\log(2)}{k-K+1}\right)^2$$
$$+ \sum_{E \in \mathfrak{D}_k(V)} \mathbb{P}_{\mathcal{F}_K(j)}\{E \in \mathcal{P}_k\} C^2 k^{-2}$$
$$\leq 16C^2 k^{-2} 2^{-\gamma(k-K)} \sum_{E \in \mathfrak{D}_k(V)} \sum_{l=1}^{2^{k-K}} l^{1-\gamma}\left(1 - \frac{\log(l)/\log(2)}{k-K+1}\right)^2$$
$$+ C^2 k^{-2} 2^{(2-\gamma)(k-K)}.$$

Since $\#\mathfrak{D}_k(V) = 4^{k-K}$, it is enough to show that

$$(3.19) \quad \sum_{l=1}^{2^{k-K}} l^{1-\gamma}\left(1 - \frac{\log(l)/\log(2)}{k-K+1}\right)^2 \leq C \frac{2^{(2-\gamma)(k-K+1)}}{(k-K+1)^2}$$

for some constant $C$ (independent of $K$). To this end note that we can compare the sum with the integral

$$\sum_{l=1}^{2^{k-K}} l^{1-\gamma}\left(1 - \frac{\log(l)/\log(2)}{k-K+1}\right)^2$$

$$\leq 2\int_1^{2^{k-K+1}} x^{1-\gamma}\left(1 - \frac{\log(x)/\log(2)}{k-K+1}\right)^2 dx$$

$$= 2(k-K+1)\log(2)\int_0^1 e^{(2-\gamma)(k-K+1)\log(2)y}(1-y)^2 dy$$

$$= 2(k-K+1)\log(2)e^{(2-\gamma)(k-K+1)\log(2)}\int_0^1 e^{(2-\gamma)(k-K+1)\log(2)y}y^2 dy$$

$$\leq \frac{4}{\log(2)^2(2-\gamma)^3}\frac{2^{(2-\gamma)(k-K+1)}}{(k-K+1)^2}. \qquad \square$$

Having proved the lemma, the proof of Proposition 3.8 is now complete.

3.4. *Lower bounds*: *proof of Proposition* 3.5.

3.4.1. *Admissible cubes, locally and globally thin points.* Fix $1 < b < a/(4-d)$ and by our choice $\gamma = (4-d)\xi/a + (4-d) - \xi$ one can choose $\varepsilon > 0$ so small that $4 - d - \gamma - 2d\varepsilon > (4-d)(b-1)\xi/a + 8b\varepsilon/(1-\varepsilon)$.

Suppose that $E \in \mathfrak{D}_k$ is an open dyadic cube and that both Brownian motions $W^1$ and $W^2$ hit the cube $E$. Then we write $\sigma_E^i := \tau^i(E, B(E, 2^{-k(1-\varepsilon)})^c)$ for the first exit times from $B(E, 2^{-k(1-\varepsilon)})$ after they first hit $E$, for $i = 1, 2$. We let

(3.20) $$S(E) = W^1([0, \sigma_E^1]) \cap W^2([0, \sigma_E^2])$$

be the intersection of the paths up to these stopping times.

DEFINITION 3.12 (*Thin points*). We define the sets

$$\mathcal{T}_{k,b} := \bigcup_{\substack{0 < r < 2^{-k} \\ 0 < \delta < 1 < \eta}} \{x \in S^* : (B(x, r\eta) \setminus B(x, r^b\delta)) \cap S = \varnothing\},$$

$$\mathcal{T}_{k,b}^{\mathrm{loc}} := \bigcup_{E \in \mathfrak{D}_k} \bigcup_{\substack{0 < r < 2^{-k} \\ 0 < \delta < 1 < \eta}} \{x \in E \cap S^* : (B(x, r\eta) \setminus B(x, r^b\delta)) \cap S(E) = \varnothing\}.$$

We call the points in $\mathcal{T}_{k,b}$ and $\mathcal{T}_{k,b}^{\mathrm{loc}}$ globally $(k,b)$-thin, respectively, locally $(k,b)$-thin.

Note that the parameters $\delta, \eta$ ensure that $\mathcal{T}_{k,b}$ and $\mathcal{T}_{k,b}^{\text{loc}}$ are open sets.
Recall that $\mathcal{P}_k$ is the set of cubes kept in the $k$th step of the percolation.

Recall from Definition 3.7 that, for $k \in \mathbb{N}$ and an open set $U \subset \mathsf{C}$,

$$\mathcal{A}_k = \{E \in \mathfrak{D}_k(U) \cap \mathcal{P}_k : \tau^i(E) < \tau^i(B_R^c),$$
$$\tau^i(E, B(E, 2^{-k(1-\varepsilon)})^c, E) > \tau^i(B_R^c), \text{ for } i = 1, 2\}$$

is the set of *admissible* cubes in $E \in \mathfrak{D}_k(U)$.

Consider the open set $U$ fixed for the moment. We subdivide $\mathfrak{D}_k(U)$ into $m_k$ disjoint subcollections $\mathfrak{D}_k(U, 1), \ldots, \mathfrak{D}_k(U, m_k)$ such that

(3.21)
$$c_1 := \inf\{2^{-2dk\varepsilon} m_k, \ k \in \mathbb{N}\} > 0,$$
$$c_2 := \sup\{2^{-2dk\varepsilon} m_k, \ k \in \mathbb{N}\} < \infty$$

and

(3.22)
$$B(E, 2^{-k(1-2\varepsilon)}) \cap F = \varnothing \quad \text{for all } E, F \in \mathfrak{D}_k(U, j), \ E \neq F,$$
$$j \in \{1, \ldots, m_k\}, \text{ and } k \in \mathbb{N}.$$

We further define

(3.23)
$$\mathcal{A}_k(j) := \mathcal{A}_k \cap \mathfrak{D}_k(U, j).$$

We now introduce a $\sigma$-algebra $\mathcal{F}_k(j)$ which makes $\#\mathcal{A}_k(j)$ measurable without using too much information of the paths inside a sufficiently large number of cubes $E \in \mathfrak{D}_k(U, j)$. The idea, as twice before, is to consider the first entrance of a path in any box $E(1)$, after that its first exit of the ball $B(E(1), 2^{-k(1-\varepsilon)})$ around $E(1)$, after that its first entrance into a new box $E(2)$ and so on. $\mathcal{F}_k(j)$ will then use information of the paths between the successive times of leaving $B(E(n), 2^{-k(1-\varepsilon)})$ and entering $B(E(n+1))$, $n \in \mathbb{N}$, as well as the information of $\mathcal{P}_k$ (the percolation at generation $k$).

We fix $j$ and for the moment suppress it in the notation. Formally for $i = 1, 2$ we introduce a sequence of (random) sets $(E^i(n) : n = 1, \ldots, m^i)$ and stopping times

$$0 = \sigma^i(0) < \rho^i(1) < \sigma^i(1) < \rho^i(2) < \cdots < \sigma^i(m^i) < \rho^i(m^i + 1)$$

by

$$\rho^i(1) := \inf\{\tau^i(E) : E \in \mathfrak{D}_k(U, j)\},$$
$$\tau^i(E^i(n)) = \rho^i(n) \quad [\text{this defines } E^i(n)],$$
(3.24) $\quad \sigma^i(n) = \tau^i(E^i(n), B(E^i(n), 2^{-k(1-\varepsilon)})^c) \quad \text{if } \rho^i(n) < \tau^i(B_R^c),$
$$\rho^i(n+1) = \inf\{\tau^i(E) : E \in \mathfrak{D}_k(U, j) \setminus \{E^i(1), \ldots, E^i(n)\}\},$$
$$m^i := \max\{n : \rho^i(n) < \tau^i(B_R^c)\}.$$

Let $\widetilde{W}_n^i(t) = W^i(t + \sigma^i(n))$, for $0 \leq t \leq \rho^i(n+1) - \sigma^i(n)$. Denote

$$\mathcal{F}_k^i(j) = \sigma(\widetilde{W}^i(t),\ t \in [0, \rho^i(n+1) - \sigma^i(n)],\ n = 0, \ldots, m^i),$$
$$\mathcal{F}_k(j) = \mathcal{F}_k^1(j) \vee \mathcal{F}_k^2(j) \vee \sigma(\mathcal{P}_k).$$

The following lemma is immediate from the construction of $\mathcal{F}_k(j)$.

LEMMA 3.13. *Admissibility of a cube is an $\mathcal{F}_k(j)$-measurable event,*

$$\{E \in \mathcal{A}_k(j)\} \in \mathcal{F}_k(j) \qquad \text{for } E \in \mathfrak{D}_k(U, j).$$

We now use Proposition 3.8 to make sure that there is a sufficiently large number of admissible cubes. Fix some $\zeta$ such that $4 - d - \gamma - 2d\varepsilon > \zeta > ((4-d)\xi/a + \varepsilon)(b-1) + 8b\varepsilon/(1-\varepsilon)$.

LEMMA 3.14.
$$\lim_{k \to \infty} \mathbb{P}\{\#\mathcal{A}_k(j) \leq 2^{\zeta k} \text{ for all } 1 \leq j \leq m_k | U \cap S^* \neq \varnothing\} = 0.$$

PROOF. Note that, for all sufficiently large $k$, the event $\#\mathcal{A}_k > c_2 2^{(\zeta + 2d\varepsilon)k}$ [recall $c_2$ from (3.21)] implies that there exists a $1 \leq j \leq m_k$ with the property that $\#\mathcal{A}_k(j) > 2^{k\zeta}$. Recalling Proposition 3.8 we obtain

$$\limsup_{k \to \infty} \mathbb{P}\{\#\mathcal{A}_k(j) \leq 2^{\zeta k} \text{ for all } 1 \leq j \leq m_k | U \cap S^* \neq \varnothing\}$$
$$\leq \limsup_{k \to \infty} \mathbb{P}\{\#\mathcal{A}_k \leq c_2 2^{(\zeta + 2d\varepsilon)k} | U \cap S^* \neq \varnothing\} = 0. \qquad \square$$

Fix $m$ such that $(1-\varepsilon)m \leq (k+5)b < (1-\varepsilon)m + 1$. For $E \in \mathfrak{D}_k(U)$ and $m$ let

$$\mathfrak{D}_{k,m}(E) = \{F \in \mathfrak{D}_m : F \subset E,\ \mathrm{dist}(F, \partial E) \geq \tfrac{5}{12} 2^{-k}\}.$$

The numbers are carefully chosen such that

(3.25) $$\bigcup_{G \in \mathfrak{D}_{k,m}(E)} G \subset B(F, \tfrac{1}{3} 2^{-k}) \subset E \qquad \text{for any } F \in \mathfrak{D}_{k,m}(E).$$

Recall the definition of $S(E)$ from (3.20).

DEFINITION 3.15 (*Successful cubes*). A cube $E$ in

$$\mathcal{S}_k(U) := \{E \in \mathcal{A}_k : \text{there exists } F \in \mathfrak{D}_{k,m}(E) \text{ such that}$$
$$F \cap S^* \neq \varnothing \text{ and } (B(F, \tfrac{1}{3} 2^{-k}) \setminus B(F, 2^{-(1-\varepsilon)m})) \cap S(E) = \varnothing\}$$

is called successful.

LEMMA 3.16. $E \cap \mathcal{T}_{k,b} \neq \varnothing$ for all $E \in \mathcal{S}_k(U)$, if $k$ is sufficiently large.

PROOF. Let $E \in \mathcal{S}_k(U)$. We first show that $E \cap \mathcal{T}_{k,b}^{\mathrm{loc}} \neq \varnothing$. Suppose $F \in \mathfrak{D}_{k,m}(E)$ satisfies the conditions in the definition of $\mathcal{S}_k(U)$. Pick $x \in F \cap S^*$, and let $r = 2^{-k-3}$, and $\eta = 2$, and $\delta = 1/2$. Then we have $B(x, r\eta) \subset B(F, \frac{1}{3}2^{-k})$ and $B(x, r^b\delta) \supset B(F, 2^{-(1-\varepsilon)m})$. Hence $x \in \mathcal{T}_{k,b}^{\mathrm{loc}}$. Finally, note that $x \in \mathcal{T}_{k,b}^{\mathrm{loc}}$ and $x \in E$ for some $E \in \mathcal{A}_k$ implies $x \in \mathcal{T}_{k,b}$. □

3.4.2. *The main step.* The following proposition is at the heart of our proof.

PROPOSITION 3.17. *Almost surely, $\mathcal{T}_{k,b}$ is dense in $S^*$ for all $k \in \mathbb{N}$.*

The rest of this section is devoted to the proof of Proposition 3.17. We have to show that, for every $k \in \mathbb{N}$ and every open set $U$ in a countable basis of the topology on $\mathsf{C}$ we have $\mathbb{P}\{\mathcal{T}_{k,b} \cap U \neq \varnothing | U \cap S^* \neq \varnothing\} = 1$. For this purpose keep $U$ fixed, as in the previous section, note that $\mathcal{T}_{k+1,b} \subset \mathcal{T}_{k,b}$, and note that it is sufficient to show that

$$\lim_{k \to \infty} \mathbb{P}\{\mathcal{S}_k(U) = \varnothing | U \cap S^* \neq \varnothing\} = 0.$$

We make the rough estimate

$$\mathbb{P}\{\mathcal{S}_k(U) = \varnothing | U \cap S^* \neq \varnothing\}$$
$$\leq \mathbb{P}\{\#\mathcal{A}_k(j) \leq 2^{\zeta k} \ \forall j = 1, \ldots, m_k | U \cap S^* \neq \varnothing\}$$
$$+ \sum_{j=1}^{m_k} \frac{\mathbb{P}\{\mathfrak{D}_k(U, j) \cap \mathcal{S}_k(U) = \varnothing | \#\mathcal{A}_k(j) > 2^{\zeta k}\}}{\mathbb{P}\{U \cap S^* \neq \varnothing\}}.$$

We know from Lemma 3.14 that the *first* term on the right-hand side converges to zero. Hence, it suffices to show that the *second* term on the right-hand side vanishes as $k \to \infty$.

For this purpose fix $k \in \mathbb{N}$, $j \in \{1, \ldots, m_k\}$ and $E \in \mathcal{A}_k(j)$. Recall that the random collection $\mathcal{A}_k(j)$ is $\mathcal{F}_k(j)$-measurable. Further recall that $\tau^i(E)$ is the time of first entry of $W^i$ into $E$ and $\tau^i(E, B(E, 2^{-k(1-\varepsilon)})^c)$ its first time to exit $B(E, 2^{-k(1-\varepsilon)})$ again. Let $\sigma_E^i := \tau^i(E, B(E, 2^{-k(1-\varepsilon)})^c) - \tau^i(E)$ and

$$V_E^i : [0, \sigma_E^i] \to B(E, 2^{-k(1-\varepsilon)}), \qquad t \mapsto W^i(t + \tau^i(E)).$$

Conditional on $\mathcal{F}_k(j)$ each $V_E^i$ is a conditioned Brownian motion with fixed start and exit points.

Write $G(E)$ for the event that there exists $F \in \mathfrak{D}_{k,m}(E)$ such that:

(a) $\dim(V_E^1[0, \sigma_E^1] \cap V_E^2[0, \sigma_E^2] \cap \Gamma[\gamma] \cap F) > 0$,
(b) $(B(F, \frac{1}{3}2^{-k}) \setminus B(F, 2^{-(1-\varepsilon)m})) \cap S(E) = \varnothing$.

By Lemma 3.6(ii), item (a) implies $F \cap S^* \neq \varnothing$, and hence $G(E)$ implies that $E \in \mathcal{S}_k(U)$. Moreover, conditional on $\mathcal{F}_k(j)$, the family $(G(E), E \in \mathcal{A}_k(j))$ is independent. Next we give a lower bound for $\mathbb{P}\{G(E)|\mathcal{F}_k(j)\}$ on $\{E \in \mathcal{A}_k(j)\}$.

LEMMA 3.18. *There exists a constant $k_0 = k_0(b,\varepsilon)$ such that, almost surely, for $k \geq k_0$,*
$$\mathbb{P}\{G(E)|\mathcal{F}_k(j)\} \geq 2^{-k(((4-d)\xi/a+\varepsilon)(b-1)+8k\varepsilon/(1-\varepsilon))} \mathbb{1}_{\{E \in \mathcal{A}_k(j)\}}.$$

PROOF. We use the notation $\mathbb{P}_{\mathcal{F}_k(j)} := \mathbb{P}\{\cdot|\mathcal{F}_k(j)\}$.

We first fix a cube $F \in \mathfrak{D}_{k,m}(E)$ and give a lower bound for the probability that $F$ satisfies the conditions (a) and (b). Note that there is a constant $C_0 > 0$ such that the event
$$H(F) := \{\tau^i(B(F, \tfrac{1}{3}2^{-k})) < \tau^i(B_R^c) \text{ for } i = 1,2,$$
$$\text{and } |W^1(\tau^1(B(F,\tfrac{1}{3}2^{-k}))) - W^2(\tau^2(B(F,\tfrac{1}{3}2^{-k})))| > \tfrac{1}{6}2^{-k}\}$$
has probability $\mathbb{P}_{\mathcal{F}_k(j)}(H(F)) > C_0$. Moreover, denote
$$I(F) := \bigcap_{i=1}^{2} \{\tau^i(B(F,\tfrac{1}{2}2^{-(1-\varepsilon)m})) < \tau^i(B_R^c)\}$$
and
$$J(F) := \bigcap_{i=1}^{2} \{\tau^i(B(F,\tfrac{1}{2}2^{-(1-\varepsilon)m}), B(E, 2^{(1-\varepsilon)k})^c)$$
$$< \tau^i(B(F,\tfrac{1}{2}2^{-(1-\varepsilon)m}), B(E, \tfrac{1}{2}2^{(1-\varepsilon)k})^c), B(F, 2^{-(1-\varepsilon)m})\}.$$

By Corollary 2.2, there exists a constant $C_1 > 0$ such that almost surely,
$$\mathbb{P}_{\mathcal{F}_k(j)}(I(F)|H(F)) > C_1 2^{2(d-2)(k-m)} \mathbb{1}_{\{E \in \mathcal{A}_k(j)\}}.$$

Now assume that $H(F)$ and $I(F)$ hold. We split each path into three pieces,
$$W^i_{(1)} : [0, \tau^i(B(F,\tfrac{1}{2}2^{-(1-\varepsilon)m})) - \tau^i(B(F,\tfrac{1}{3}2^{-k}))] \to \mathbb{R}^d,$$
$$W^i_{(1)}(t) = W^i(\tau^i(B(F,\tfrac{1}{3}2^{-k})) + t),$$
$$W^i_{(2)} : [0, \tau^i(B(F,\tfrac{1}{2}2^{-(1-\varepsilon)m}), \partial B(F, 2^{-(1-\varepsilon)m})) - \tau^i(B(F,\tfrac{1}{2}2^{-(1-\varepsilon)m}))] \to \mathbb{R}^d,$$
$$W^i_{(2)}(t) = W^i(\tau^i(B(F,\tfrac{1}{2}2^{-(1-\varepsilon)m})) + t),$$
$$W^i_{(3)} : [0, \tau^i(B(F,\tfrac{1}{2}2^{-(1-\varepsilon)m}), \partial B(E, \tfrac{1}{2}2^{-k(1-\varepsilon)}))$$
$$- \tau^i(B(F,\tfrac{1}{2}2^{-(1-\varepsilon)m}), \partial B(F, 2^{-(1-\varepsilon)m}))] \to \mathbb{R}^d,$$
$$W^i_{(3)}(t) = W^i(\tau^i(B(F,\tfrac{1}{2}2^{-(1-\varepsilon)m}), \partial B(F, 2^{-(1-\varepsilon)m})) + t).$$

We now form the packet consisting of the first and last part of the $i$th motion $(i = 1, 2)$,

$$\mathcal{W}^i_{1\cup 3} = W^i_{(1)}[0, \tau^i(B(F, \tfrac{1}{3}2^{-(1-\varepsilon)m})) - \tau^i(B(F, \tfrac{1}{3}2^{-k}))]$$
$$\cup W^i_{(3)}[0, \tau^i(B(F, \tfrac{1}{2}2^{-(1-\varepsilon)m}), \partial B(E, \tfrac{1}{2}2^{-k(1-\varepsilon)}))$$
$$- \tau^i(B(F, \tfrac{1}{2}2^{-(1-\varepsilon)m}), \partial B(F, 2^{-(1-\varepsilon)m}))].$$

Using Lemma 2.9 and the subsequent Remark 2.10 we get that for a suitable constant $k_0$, depending only on $b$ and $\varepsilon$, and for all $k \geq k_0$,

$$(3.26) \quad \mathbb{P}_{\mathcal{F}_k(j)}\{\mathcal{W}^1_{1\cup 3} \cap \mathcal{W}^2_{1\cup 3} = \varnothing | H(F) \cap I(F)\} \geq 2^{(\xi+\varepsilon)(k-m)} \mathbb{1}_{\{E \in \mathcal{A}_k(j)\}}.$$

Observe that this event as well as $H(F)$ and $I(F)$ are measurable with respect to the $\sigma$-field $\mathcal{G} := \sigma(W^1_{(1)}, W^2_{(1)}, W^1_{(3)}, W^2_{(3)})$. Now define the ranges of the middle pieces

$$\mathcal{W}^i_2 := W^i_{(2)}[0, \tau^i(B(F, \tfrac{1}{2}2^{-(1-\varepsilon)m}), \partial B(F, 2^{-(1-\varepsilon)m})) - \tau^i(B(F, \tfrac{1}{2}2^{-(1-\varepsilon)m}))]$$
$$\text{for } i = 1, 2,$$

where we agree that $W^i_{(1)} = W^i$ and $W^i_{(3)} = W^i_{(2)} = \varnothing$ if $H(F)$ does not occur.

We can find a constant $C_2 > 0$ such that almost surely,

$$\mathbb{P}_{\mathcal{F}_k(j)}\{\dim(\mathcal{W}^1_2 \cap \mathcal{W}^2_2 \cap F \cap \Gamma[\gamma]) > 0 | \mathcal{G}\} \geq C_2 2^{\gamma(k-m)-2m\varepsilon} \mathbb{1}_{\{E \in \mathcal{A}_k(j)\}} \mathbb{1}_{H(F) \cap I(F)}.$$

Indeed, on $\{E \in \mathcal{A}_k(j)\}$ the probability that $F$ is retained in the percolation equals $2^{\gamma(k-m)}$ and on $I(F)$ the probability that both motions hit the cube $F$ is no smaller than a constant multiple of $2^{-2m\varepsilon}$. Given that both motions hit $F$ and $F \in \mathcal{P}_m$, by Lemma 3.3(iii) there is a positive probability that $\mathcal{W}^1_2 \cap \mathcal{W}^2_2 \cap F \cap \Gamma[\gamma]$ has positive dimension, and it is easy to see by Brownian scaling that this probability does not depend on the scale, that is, on $m$.

Finally, a simple argument shows that, for any $x \in \partial B(E, \tfrac{1}{2}2^{-(1-\varepsilon)k})$,

$$\mathbb{P}_{\mathcal{F}_j(k)}(J(F)|\mathcal{G}) \geq c^2 \mathbb{P}_x\{\tau(B(E, 2^{-(1-\varepsilon)k})^c) < \tau(B(E, 2^{-k}))\}^2 \mathbb{1}_{I(F)},$$

where [recall that $\mathcal{U}_r$ is the uniform distribution on $\partial B(0, r)$]

$$c = \inf_{x \in \partial B(0,1)} \inf_{y \in \partial B(0,2)} \frac{\mathbb{P}_x\{W(\tau(B(0,2))) \in dy\}}{\mathcal{U}_2(dy)} > 0.$$

Hence, by Lemma 2.1 there exists a constant $C_3 = C_3(\varepsilon) > 0$ such that

$$\mathbb{P}_{\mathcal{F}_k(j)}(J(F)|\mathcal{G}) \geq C_3 \mathbb{1}_{I(F)}.$$

Note that $\{\mathcal{W}^1_{1\cup 3} \cap \mathcal{W}^2_{1\cup 3} = \varnothing\} \cap J(F)$ implies item (b).

Altogether, the probability that a *fixed* cube $F \in \mathfrak{D}_{k,m}(E)$ satisfies (a) and (b) is (with $C = C_1 \cdot C_2 \cdot C_3$)

$$\mathbb{P}_{\mathcal{F}_k(j)}\{F \text{ satisfies (a) and (b)}\} \tag{3.27}$$
$$\geq C 2^{(\gamma+2(d-2))(k-m)+(\xi+\varepsilon)(k-m)-2m\varepsilon} \mathbb{1}_{\{E \in \mathcal{A}_k(j)\}},$$

almost surely. Note now that if $F$ satisfies (b) then, by (3.25) no cube in $\mathfrak{D}_{k,m}(E)$ which does not intersect $B(F, 2^{-(1-\varepsilon)m})$ satisfies (a). Hence (a) and (b) are satisfied by at most a constant multiple of $2^{dm\varepsilon}$ cubes in $\mathfrak{D}_{k,m}(E)$ simultaneously. As the total number of cubes $F \in \mathfrak{D}_{k,m}(E)$ is at least a constant multiple of $2^{d(m-k)}$, we have constants $C_4, C_5 > 0$ such that, on $\{E \in \mathcal{A}_k(j)\}$,

$$\mathbb{P}_{\mathcal{F}_k(j)}(G(E)) \geq C_4 2^{-dm\varepsilon} \sum_{F \in \mathfrak{D}_{k,m}(E)} \mathbb{P}_{\mathcal{F}_k(j)}\{F \text{ satisfies (a) and (b)}\} \mathbb{1}_{\{E \in \mathcal{A}_k(j)\}}$$

$$\geq C_5 2^{(\gamma+d-4)(k-m)+(\xi+\varepsilon)(k-m)-5m\varepsilon} \mathbb{1}_{\{E \in \mathcal{A}_k(j)\}}.$$

By definition of $\gamma$, we have $\gamma + d - 4 + \xi = (4-d)\xi/a < 3$ and, by definition of $m$, $k - m \geq k(1-b) - m\varepsilon + 5b$, and $m \leq (k+5)b/(1-\varepsilon)$, which gives the claimed lower bound. □

Using this lemma and the conditional independence of the family $(G(E), E \in \mathcal{A}_k(j))$ we get for $k \geq k_0$, almost surely,

$$\mathbb{P}\{\mathfrak{D}_k(U, j) \cap \mathcal{S}_k(U) = \varnothing | \mathcal{F}_k(j)\} \leq (1 - 2^{-k(((4-d)\xi/a+\varepsilon)(b-1)+8b\varepsilon/(1-\varepsilon))})^{\#\mathcal{A}_k(j)}.$$

Finally, this gives for $k \geq k_0$

$$\sum_{j=1}^{m_k} \mathbb{P}\{\mathfrak{D}_k(U, j) \cap \mathcal{S}_k(U) = \varnothing | \#\mathcal{A}_k(j) > 2^{\zeta k}\}$$

$$\leq c_2 2^{2dk\varepsilon}(1 - 2^{-k(((4-d)\xi/a+\varepsilon)(b-1)+8b\varepsilon/(1-\varepsilon))})^{2^{\zeta k}}.$$

By our choice of parameters $\zeta > ((4-d)\xi/a + \varepsilon)(b-1) + 8b\varepsilon/(1-\varepsilon)$ and hence the right-hand side converges to 0 as $k \to \infty$. This completes the proof of Proposition 3.17.

3.4.3. *Completing the proof of Proposition* 3.5: *a density argument.* Recall that $b \in (1, a/(4-d))$ was chosen arbitrarily. Note that

$$\bigcap_{k=1}^{\infty} \mathcal{T}_{k,b} \subset \mathcal{T}(b(4-d), R) \cap S^*.$$

Indeed, if $x \in \bigcap_{k=1}^{\infty} \mathcal{T}_{k,b}$, then $x \in S^*$ and there exists a sequence $r_k \downarrow 0$ with

$$(B(x, r_k) \setminus B(x, r_k^b)) \cap S = \varnothing.$$

Then, by [13], (1.17), if $d = 3$ and by [4], (1.6), if $d = 2$, there exists a constant $C > 0$ such that for sufficiently large $k \in \mathbb{N}$,

$$\ell(B(x, r_k)) = \ell(B(x, r_k^b)) \leq C r_k^{b(4-d)} [\log(1/r_k)]^2,$$

and hence $x \in \mathcal{T}(b(4-d), R)$.

Clearly $\mathcal{T}_{k,b} \subset \mathcal{T}_{k,c}$ for $b > c$, and $\mathcal{T}(a, R) = \bigcap_{b < a} \mathcal{T}(b, R)$. Thus

$$\bigcap_{n=1}^{\infty} \bigcap_{k=1}^{\infty} \mathcal{T}_{k, a/(4-d) - (1/n)} \subset \mathcal{T}(a, R) \cap S^*.$$

Next recall that $\mathcal{T}_{k, a/(4-d) - (1/n)}$ is relatively open in $S^*$ and, by Proposition 3.17, also dense in $S^*$ for any $k, n$. As $S^*$ is compact, hence complete, one can infer from Baire's theorem that $\bigcap_{k,n \in \mathbb{N}}^{\infty} \mathcal{T}_{k, a/(4-d) - (1/n)}$ is dense in $S^*$ almost surely. Hence $\mathbb{P}\{\mathcal{T}(a, R) \cap S^* \neq \varnothing | S^* \neq \varnothing\} = 1$ and, since $S^* \subset \Gamma[\gamma] \cap S \subset \Gamma[\gamma]$, we have

$$\begin{aligned}
\mathbb{P}\{\mathcal{T}(a, R) \cap \Gamma[\gamma] \neq \varnothing\} &\geq \mathbb{P}\{\mathcal{T}(a, R) \cap S^* \neq \varnothing\} \\
&= \mathbb{P}\{\mathcal{T}(a, R) \cap S^* \neq \varnothing | S^* \neq \varnothing\} \mathbb{P}\{S^* \neq \varnothing\} \\
&= \mathbb{P}\{S^* \neq \varnothing\} > 0.
\end{aligned}$$

This completes the proof of Proposition 3.5.

**Acknowledgments.** We thank the *Mathematisches Forschungsinstitut Oberwolfach* and the organizers of the Miniworkshop "Stochastische Prozesse in zufälligen Medien" in May 2002 where our joint work was initiated.

MATHEMATISCHES INSTITUT
JOHANNES GUTENBERG-UNIVERSITÄT MAINZ
STAUDINGERWEG 9
55099 MAINZ
GERMANY
E-MAIL: [math@aklenke.de](math@aklenke.de)

DEPARTMENT OF MATHEMATICAL SCIENCES
UNIVERSITY OF BATH
CLAVERTON DOWN
BATH BA2 7AY
UNITED KINGDOM
E-MAIL: [maspm@bath.ac.uk](maspm@bath.ac.uk)